\documentclass{amsart}
\makeatletter
\def\l@subsection{\@tocline{2}{0pt}{2pc}{5pc}{}}
\makeatother

\usepackage{bm}
\usepackage{amsmath,amssymb,amsthm,mathtools}
\usepackage{mathrsfs}
\usepackage{booktabs,array}
\usepackage{enumitem}
\usepackage{needspace}
\usepackage{tikz-cd}
\usetikzlibrary{calc, decorations.pathreplacing}
\usepackage{hyperref}
\usepackage{comment}

\hypersetup{
	colorlinks=true,
	linkcolor=black,
	citecolor=blue,
	urlcolor=blue
}

\theoremstyle{plain}
\newtheorem{theorem}{Theorem}[section]
\newtheorem{proposition}[theorem]{Proposition}
\newtheorem{lemma}[theorem]{Lemma}
\newtheorem{corollary}[theorem]{Corollary}

\theoremstyle{definition}
\newtheorem{definition}[theorem]{Definition}
\newtheorem{example}[theorem]{Example}

\theoremstyle{remark}

\newcommand{\ZZ}{\mathbb{Z}}
\newcommand{\QQ}{\mathbb{Q}}
\newcommand{\RR}{\mathbb{R}}
\newcommand{\CC}{\mathbb{C}}
\newcommand{\GG}{\mathbb{G}}
\newcommand{\PP}{\mathbb{P}}

\newcommand{\cM}{\mathcal{M}}
\newcommand{\cX}{\mathcal{X}}
\newcommand{\cS}{\mathcal{S}}
\newcommand{\cT}{\mathcal{T}}
\newcommand{\cO}{\mathcal{O}}
\newcommand{\cV}{\mathcal{V}}
\newcommand{\cD}{\mathcal{D}}

\newcommand{\Sec}{\operatorname{Sec}}
\newcommand{\Span}{\operatorname{Span}}

\newcommand{\rk}{\operatorname{rk}}

\newcommand{\Vol}{\operatorname{Vol}}

\newcommand{\Hom}{\operatorname{Hom}}

\newcommand{\deq}{\coloneqq}

\newcommand{\STP}{\mathrm{STP}}
\newcommand{\tSTP}{\widetilde{\mathrm{STP}}}
\newcommand{\cSTP}{\mathcal{STP}}
\newcommand\bkp{\bar{\kappa}}
\newcommand\tbkp{\widetilde{\kappa}}
\newcommand\bkpt{\kappa^{\mathrm{tor}}}
\newcommand\kp{\kappa}

\title[Kappa Classes on the KSBA Space of Stable Toric Pairs]{Kappa Classes on the KSBA Space \\of Stable Toric Pairs}
\date{}
\author{Byeol Han}
\email{Byeol@uga.edu}
\address{Department of Mathematics, University of Georgia, Athens, GA 30602, USA}
\begin{document}
	
	\begin{abstract}
		We study kappa classes on the normalization of the main
		component of the Koll\'ar--Shepherd-Barron--Alexeev
		(KSBA) moduli space of stable toric pairs. Using a toric model obtained by a lattice modification, we compute all kappa classes on the moduli space by giving a combinatorial formula for their intersection numbers with all torus-invariant cycles. In particular, we identify the first normalized kappa class with the GKZ polarization of the toric variety corresponding to the secondary polytope. For interval configurations, we also obtain a closed formula for all top intersection products of normalized kappa classes.
	\end{abstract}
	\maketitle
	
	\tableofcontents

	\section{Introduction}
	
	The kappa classes, also known as the Miller--Morita--Mumford (MMM)
	classes, are fundamental tautological classes on the moduli spaces
	$\overline{M}_{g,n}$ of stable curves. They were introduced by Mumford
	in \cite{Mum} and, independently from a topological perspective, by
	Morita \cite{Mor} and were further studied by Miller \cite{Mil}.
	Their importance is reflected in the topology of the moduli space.
	For example, the stable rational cohomology ring of $M_g$ is
	$\QQ[\kp_1,\kp_2,\ldots]$ by Mumford's conjecture, proved by Madsen
	and Weiss \cite{MW07}. 
	
	Donaldson \cite{Don23} proposed an approach to extending
	this theory to KSBA moduli spaces of surfaces, which generalize the
	moduli spaces of stable curves to higher dimensions. This approach is being further developed by Jiang in
	\cite{JiangVFC,JiangEnumerativeKSBA}, while Alexeev
	\cite{Ale23,AlexeevKappaBasic} developed kappa classes on
	general KSBA moduli spaces. For a universal KSBA family
	$f:(\mathcal X,\mathcal D)\to\mathcal M$ of relative dimension~$d$,
	these classes are defined by
	$$
	\kp_r
	=
	f_*\left(
	c_1\bigl(
	\mathcal O_{\mathcal X}
	(K_{\mathcal X/\mathcal M}+\mathcal D)
	\bigr)^{d+r}
	\right)
	\in A^r(\mathcal M)_\QQ,
	$$
	where $f_*$ is the Gysin pushforward.
	
We compute these classes for stable toric pairs.
	Let $A\subset M\cong\ZZ^d$ be a full-dimensional finite
	lattice configuration, $Q = \operatorname{Conv}(A)$, and let $\cSTP$ be the normalization of the main component of the moduli stack of stable toric pairs associated to~$A$.
	Its coarse moduli space $\STP$ is the toric variety
	associated with $\Sec(A)$ \cite{Alexeev02}.
	The corresponding KSBA stable pairs $(X,\Delta+\varepsilon B)$,
	where $\Delta$ is the toric boundary and
	$0<\varepsilon\ll1$, define classes $\kp_r(\varepsilon)$ on $\cSTP$.
	We write $\kp_r=\varepsilon^{-(d+r)}\kp_r(\varepsilon)$
	for the normalized class and
	$\bkp_r\in A^r(\STP)_\QQ$
	for its descent to the coarse moduli space. For a logarithmic approach to the moduli of stable toric varieties, see Ascher and Molcho~\cite{AscherMolcho16}.
	
	Our computation uses an auxiliary toric equidimensional family over $\STP$ obtained from a fiber-polytope construction, which we call the effective
	toric model. Over the orbit of $\STP$ corresponding to a regular subdivision
	$\cS$ of $(Q, A)$, the reduced irreducible components of the fiber
	are indexed by the maximal cells of $\cS$.
	The effective family does not need to be flat, and its fibers
	need not be reduced. Its generic fiber has character lattice
	$\Lambda_A=\ZZ(A-A)$, which may be a proper sublattice of $M$. Thus we modify the lattices to obtain a flat family
	with reduced fibers and fiber character lattice $M$
	over a finite toric cover of $\STP$. This provides a correct family
	of stable toric pairs on which kappa classes can be
	computed by toric intersection theory.

	The intersection-theoretic framework is provided by
	operational Chow theory \cite{FM}.
	Vistoli identifies the rational
	operational Chow cohomology of the moduli stack with
	that of its coarse moduli space in
	\cite[Proposition~6.1]{Vis}.
	In the KSBA setting, Alexeev shows that the intrinsic
	kappa classes descend uniquely to the coarse moduli space
	\cite[Proposition~3.10]{AlexeevKappaBasic}.
	Since $\STP$ is a complete toric variety, its operational
	Chow classes are determined by their intersection numbers
	with torus-invariant orbit closures \cite{FS97}.
	We therefore determine the normalized kappa classes
	by the following intersection formula.
	
	\Needspace{9\baselineskip}
	\begin{theorem}
		Let $A\subset M\cong\ZZ^d$ be a full-dimensional finite
		lattice configuration, and let $r\geq 0$. For every
		$r$-dimensional torus-invariant orbit closure
		$V_{\cS}\subset\STP$ corresponding to a regular subdivision
		$\cS$ of $A$, we have
		$$
		\deg\left(
		\bkp_r\cap[V_{\cS}]
		\right)
		=
		\sum_{\substack{C\in\cS^{\max}\\|C|=d+r+1}}
		g_C,
		$$
		where 
		$g_C=[M:\ZZ(C-C)].$
	\end{theorem}
    Our main result gives a formula for the degrees of all normalized kappa classes on torus-invariant orbit closures. The formula applies even when $\Lambda_A \neq M.$ The corresponding intersection numbers on the moduli stack then follow from coarse descent.
	\begin{corollary}
		Let $\cV_{\cS}\subset\cSTP$ be the reduced integral substack lying over
		$V_{\cS}$, and let $e_{\cS}$ be the order of its generic stabilizer. Then 
		$$
		\deg\left(
		\kp_r\cap[\cV_{\cS}]
		\right)
		=
		\frac{1}{e_{\cS}}
		\sum_{\substack{C\in\cS^{\max}\\|C|=d+r+1}}
		g_C.
		$$
	\end{corollary}
	
	The case $r=1$ has a further interpretation in terms of the geometry
	of the secondary polytope. We identify the first normalized kappa class with the Gel'fand--Kapranov--Zelevinsky polarization in \cite{GKZ}:
	
	\begin{theorem}
		$$
		\bkp_1
		=
		c_1(L_{\mathrm{GKZ}})
		\quad\text{in }A^1(\STP)_\QQ,
		$$
		where $L_{\mathrm{GKZ}}$ is the ample line bundle associated with $\Sec(A)$.
	\end{theorem}
	As a complementary computation, we work directly with the
	effective toric model. Since this family is not flat in general,
	the usual construction of kappa classes using the flat
	pushforward of operational Chow classes does not apply
	directly. Instead, we use the tropical pushforward of
	weighted fans given in \cite{GKM} to define naive toric classes
	$\bkpt_r\in A^r(\STP)_\QQ$.
	The two constructions are related by the following identity.
	\begin{theorem}
		$g_A\,\bkpt_r=\bkp_r$, where $g_A=[M:\ZZ(A-A)]$.
	\end{theorem}
Finally, for the interval configuration
$A_N=\{0,\ldots,N\}\subset\ZZ$, with $N\geq3$,
write $\STP_N=X_{\Sec(A_N)}$.
We obtain a closed formula for all top intersection products of normalized kappa classes.
		\begin{theorem}
		Let $r_1,\ldots,r_s$ be positive integers satisfying
		$r_1+\cdots+r_s=N-1.$ Then
		$$
		\int_{\STP_N}
		\prod_{i=1}^{s}\bkp_{r_i}
		=
		N^{s-2}
		\prod_{j=1}^{s}(r_j+1).
		$$
	\end{theorem}
	
	The paper is organized as follows. In
	Section~\ref{sec:stable-toric-pairs}, we recall secondary polytopes
	and the moduli theory of stable toric pairs and construct the
	associated KSBA classifying morphism. We also recall the definition of kappa classes on the KSBA moduli space. In Section~\ref{sec:toric-models}, we construct the
	effective toric model and its Kummer modification, and we define the corresponding kappa classes on the Kummer model. In
	Section~\ref{sec:intrinsic-kummer}, we use the Kummer family to prove the main intersection
	formula, and in Section~\ref{sec:alternative-toric}, we compute the auxiliary toric kappa classes directly on the effective model and compare them with the intrinsic classes.
	In Section~\ref{sec:GKZ}, we identify the first normalized kappa class
	with the GKZ polarization. Finally, in Section~\ref{sec:interval}, we apply the main result to derive explicit formulas for the interval configuration. Throughout this paper, we work over the field $k = \CC.$
	
	\vspace{\baselineskip}
	\noindent\textbf{Acknowledgements.} The author is grateful to her advisor, Valery Alexeev, for his insightful comments and helpful discussions.
	
	\vspace{\baselineskip}
	\noindent\textbf{AI Disclosure.} 
	During the development of this paper, the author used OpenAI's
	ChatGPT, primarily with GPT-5.6 Sol and GPT-6 Astra, as an
	interactive research aid, particularly for exploring the explicit
	combinatorial calculations in Section~\ref{sec:interval}.

	\section{Stable Toric Pairs and Kappa Classes}\label{sec:stable-toric-pairs}
	The moduli space of stable toric pairs was constructed by Alexeev in \cite{Alexeev02} and extended to stable spherical varieties in \cite{AlexeevBrion06}. We recall the moduli problem using the concrete toric terminology of \cite{Weighted}.
	
	\subsection{Stable Toric Pairs}
	Let $M \cong \ZZ^d$ be a lattice, let $T$ be the algebraic torus with character lattice $M$, and let $A \subset M$ be a full-dimensional finite lattice configuration. Set $Q = \operatorname{Conv}(A) \subset M_{\RR}.$ The basic object in the moduli problem is the polarized toric variety associated with $Q$, together with the sections indexed by ~$A$. Degenerations of such toric varieties are encoded by subdivisions of $Q$ together with the subsets of $A$ occurring on their components.

	Throughout, we work in the framed category of \cite{Alexeev02}: the torus $T$ and its character
	lattice $M$ are fixed, and isomorphisms are
	$T$-equivariant for this fixed torus action.
   We also keep the labels $a\in A$ fixed, so isomorphisms
preserve the indexing of the weights by $A$.

	\begin{definition}
		A {\it subdivision of $(Q, A)$} is a family $\cS = \{(Q_i, A_i): i \in I\}$ such that
		\begin{enumerate}
			\item each $A_i$ is a full-dimensional subset of $A$ and $Q_i = \operatorname{Conv} A_i$;
			\item any intersection $Q_i \cap Q_j$ is a face of both $Q_i$ and $Q_j,$ and
			$$A_i \cap (Q_i \cap Q_j) = A_j \cap (Q_i \cap Q_j);$$
			\item $Q = \cup_{i \in I} Q_i.$
		\end{enumerate}
	\end{definition}

Each pair $(Q_i,A_i)$ is called a \emph{marked cell},
and $A_i$ is its set of marked points.
The set $A_i$ is part of the data and need not equal
$A\cap Q_i$. We record a subdivision by its maximal marked cells and write $\cS^{\max}$ for its maximal cells.
We often identify a marked cell $(Q_i,A_i)$ with its set of marked points $A_i$.
	\begin{definition}
		A {\it polarized stable toric variety} is a pair $(X,L),$ where $X$ is a projective seminormal variety and $L$ is a $T$-linearized ample line bundle such that $(X_i, L_i = L|_{X_i})$ is a polarized toric variety for every irreducible component $X_i$ of ~$X$.
	\end{definition}
	
	Every irreducible component $(X_i, L_i)$ corresponds to a lattice polytope, say $Q_i$. The intersection $X_i\cap X_j$ is a closed $T$-invariant subset of both components and is therefore a union of $T$-orbits. Combinatorially, it corresponds to a union of faces of $Q_i$ and $Q_j$.
	
	\begin{definition}
		A \emph{stable toric variety over $\PP^N$}
		is a polarized stable toric variety $(X,L)$, together with a finite $T$-equivariant morphism
		$
		p:X\to\PP^N
		$
		such that
		$
		L\cong p^*\cO_{\PP^N}(1)
		$
		as $T$-linearized line bundles.
	\end{definition}
	
	Given a polarized stable toric variety $(X,L)$, let
	$s\in H^0(X,L)$ be a section that does not vanish identically on
	any irreducible component, and write $s$ as a sum of
	$T$-eigensections. Then $B=(s=0)$ is an effective Cartier
	divisor, and the eigensections define a rational map to projective
	space. By \cite{Weighted}, this map is regular if and only if
	$B$ contains no $T$-orbit.
	\begin{definition}
		A {\it stable toric pair} $(X,B)$ is a polarized stable toric variety $(X, \cO(B))$, where $B$ is an effective ample Cartier divisor containing no \mbox{$T$-orbit of $X$}.
	\end{definition}
	
	For the fixed configuration $A$, let $\PP^A=\PP(k^A)$
	with the diagonal $T$-action whose coordinate sections have
	weights indexed by $A$.
	
	A stable toric pair $(X,B)$, with $B=(s=0)$ and the weights occurring in $s$ contained in $A$, determines a finite $T$-equivariant morphism
	$
	p_s:X\longrightarrow\PP^A
	$
	by the eigensections of $s$. Conversely, if
	$p:X\to\PP^A$ is a stable toric variety over $\PP^A$, then
	$B=p^{-1}(H_A)$, where $H_A=\left\{\sum_{a\in A}z_a=0\right\}$,
	defines the corresponding stable toric pair. Hence the two moduli
	problems are equivalent.

	\begin{definition}
	A \emph{stable toric pair associated to $A$} is a stable
	toric pair for which $\cS=\{(Q_i,A_i)\}$ is a subdivision
	of $(Q,A)$, with the components $X_i$ glued along the
	torus-invariant subvarieties corresponding to common
	faces of $\cS$.
	\end{definition}
	
	By the description above, every stable toric pair associated to $A$ comes equipped with a finite morphism to $\PP^A$; see \cite{Alexeev02}. This morphism is the starting point for the
	fiber-polytope construction of an equidimensional
	family in Section~\ref{sec:toric-models}.

	\subsection{GKZ Secondary Polytopes}
	We now recall the secondary polytope, whose faces correspond to the regular subdivisions of $(Q, A)$. Most of the material in this subsection is based on \cite[Chapter 7]{GKZ}.
	
	\begin{definition}
		A subdivision $\cS=\{(Q_i,A_i):i\in I\}$ of $(Q,A)$ is called {\it regular} if there exists a height function $h:A\to\RR$ such that each $Q_i$ is obtained by projecting a lower face $F_i$ of the convex hull of the lifted points
		$$
		\{(a,h(a))\mid a\in A\}\subset\RR^{d+1}
		$$
		to $\RR^d$, and
		$$
		A_i=\{a\in A:(a,h(a))\in F_i\}.
		$$
	\end{definition}
	
	Let $\cS=\{(Q_i,A_i)\}$ and $\cS'=\{(Q_j',A_j')\}$ be two subdivisions of $(Q,A).$ We say that $\cS$ {\it refines} $\cS'$ if the collection of $(Q_i,A_i)$ such that $Q_i\subset Q_j'$ forms a subdivision of $(Q_j',A_j')$ for every $j.$ In this sense, a triangulation is a maximal refinement of a subdivision.
	
	\begin{definition}
		A {\it triangulation of $(Q,A)$} is a triangulation of $Q$ into simplices with vertices in $A$.
	\end{definition}
Throughout the paper, $\Vol_{\ZZ}$ denotes normalized
$d$-dimensional volume with respect to $M\cong\ZZ^d$.
Thus a unimodular $d$-simplex has volume~$1$, and a
lower-dimensional polytope has volume~$0$.
	\begin{definition}
		The {\it GKZ vector of $\cT$} is the integral vector
		$(\varphi_\cT(a))_{a\in A}\in\ZZ^A$, where
		$$
		\varphi_\cT(a)
		=
		\sum_{\substack{\sigma\in\cT^{\max}\\
				a\in\operatorname{Vert}(\sigma)}}
		\Vol_{\ZZ}(\sigma).
		$$
		
		The {\it secondary polytope $\Sec(A)$ of $A$} is the convex hull in $\RR^A$ of the GKZ vectors $\varphi_\cT$ for all regular triangulations $\cT$ of $(Q,A)$. The polarization on the toric variety associated with $\Sec(A)$ is called the GKZ line bundle and is denoted by $L_{\mathrm{GKZ}}.$
	\end{definition}
	
	With the lower-hull convention, if a height function
	$h:A\to\RR$ induces a regular subdivision $\cS$, then the
	face of $\Sec(A)$ on which the linear functional
	\mbox{$x\longmapsto \langle h,x\rangle$}
	is minimized is precisely the face corresponding to $\cS$.
	In particular, if $h$ induces a regular triangulation $\cT$,
	then $\varphi_{\cT}$ is a minimizing vertex.
	
	\begin{proposition}
		\begin{enumerate}
			\item The secondary polytope $\Sec(A)$ has dimension $|A|-d-1.$
			\item There is a one-to-one correspondence between the vertices of $\Sec(A)$ and the regular triangulations of $A$.
			\item The faces of $\Sec(A)$ are the convex hulls of the GKZ vectors $\varphi_\cT$ for triangulations $\cT$ refining a fixed regular subdivision $\cS$.
		\end{enumerate}
	\end{proposition}
	
	\begin{definition}
		A subset $Z\subset A$ is called a {\it circuit} if it is minimally affinely dependent, that is, $Z$ is affinely dependent while every proper subset of $Z$ is affinely independent.
	\end{definition}
	
	The relation between degenerations of stable toric pairs and subdivisions described above is reflected in the structure of the moduli space.
	
	\begin{theorem}[\cite{Alexeev02,AlexeevBrion06}]
		Let $A\subset M\cong\ZZ^d$ be a full-dimensional finite lattice configuration.
		There exists a proper algebraic stack $\cSTP^{\mathrm{gen}}$ with finite stabilizers parametrizing stable toric pairs associated to $A$. The stack $\cSTP^{\mathrm{gen}}$ admits a projective coarse moduli space $\STP^{\mathrm{gen}}$. Moreover, the normalization of the main irreducible component of $(\STP^{\mathrm{gen}})_{\mathrm{red}}$ is isomorphic to the toric variety associated with $\Sec(A)$.
	\end{theorem}
	
	Here the main component is the closure of the locus of
	normal stable toric pairs associated to $A$. From now on, we work with the normalizations of the main components and write
	$$
	\cSTP\deq(\cSTP^{\mathrm{gen}}_{\mathrm{main}})^\nu
	\quad\text{and}\quad
	\STP\deq(\STP^{\mathrm{gen}}_{\mathrm{main}})^\nu.
	$$

\subsection{Comparison with the KSBA Moduli Space}
	In this subsection, we compare the normalized main component of the moduli space
	of stable toric pairs associated to $A$ with the corresponding locus in the
	KSBA moduli space. Recall that a KSBA stable pair is a
	semi-log-canonical (slc) pair $(X,D)$ such that $K_X+D$ is
	$\QQ$-Cartier and ample.

\begin{proposition}\label{prop:relative-log-canonical}
	Let
	$
	f:(\cX_A,\mathcal B_A)\longrightarrow\cSTP
	$
	be the universal family of stable toric pairs associated to $A$,
	and let $\Delta_{\cX_A/\cSTP}$ denote the relative toric boundary.
	Then
	$$
	K_{\cX_A/\cSTP}+\Delta_{\cX_A/\cSTP}
	\sim_{\QQ}0.
	$$
\end{proposition}

\begin{proof}
	Consider the Kummer family
	$\widetilde\pi:\widetilde X\to\tSTP$
	constructed in Section~\ref{sec:toric-models}.
	Its classifying morphism gives a base change diagram
	$$
	\begin{tikzcd}
		\widetilde X
		\arrow[r,"\bar\chi"]
		\arrow[d,"\widetilde\pi"']
		&
		\cX_A
		\arrow[d,"f"]
		\\
		\tSTP
		\arrow[r,"\chi"']
		&
		\cSTP.
	\end{tikzcd}
	$$
	
Write
$
\Delta_{\widetilde X}^{\mathrm{tor}}
=
\widetilde\Delta_{\mathrm{hor}}
+
\widetilde\Delta_{\mathrm{ver}}
$
for the horizontal and vertical parts with respect to~$\widetilde\pi$.
Here $\widetilde\Delta_{\mathrm{hor}}$ is the relative
toric boundary.
Since $\widetilde\pi$ is flat with reduced fibers,
$$
\widetilde\Delta_{\mathrm{ver}}
=
\widetilde\pi^*\Delta_{\tSTP}^{\mathrm{tor}}.
$$
For a toric variety, the log canonical divisor is trivial, and so 
$$
\begin{aligned}
	K_{\widetilde X/\tSTP}
	+\widetilde\Delta_{\mathrm{hor}}
	&=
	K_{\widetilde X}
	-\widetilde\pi^*K_{\tSTP}
	+\Delta_{\widetilde X}^{\mathrm{tor}}
	-\widetilde\Delta_{\mathrm{ver}}\\
	&=
	(K_{\widetilde X}+\Delta_{\widetilde X}^{\mathrm{tor}})
	-\widetilde\pi^*
	(K_{\tSTP}+\Delta_{\tSTP}^{\mathrm{tor}})\\
	&=0.
\end{aligned}
$$
As shown in the proof of
Proposition~\ref{prop:kummer-kappa-pullback},
$\rho\circ\chi=\mu$ is finite and surjective.
Since $\rho$ is the coarse moduli morphism and has
finite diagonal, $\chi$ and its base change $\bar\chi$
are finite and surjective.

Since $f$ is flat with reduced fibers, by the Riemann--Hurwitz formula, we have
$$
\bar\chi^*(K_{\cX_A/\cSTP}+\Delta_{\cX_A/\cSTP})
=
K_{\widetilde X/\tSTP}
+\widetilde\Delta_{\mathrm{hor}}
\sim0.
$$
Thus $(\deg\bar\chi)(K_{\cX_A/\cSTP}+\Delta_{\cX_A/\cSTP})\sim0,$  and so $K_{\cX_A/\cSTP}+\Delta_{\cX_A/\cSTP} \sim_{\QQ}0$.
\end{proof}

	\begin{corollary}\label{cor:ksba-family}
		For every sufficiently small rational
		$\varepsilon > 0 $, the family
		$$
		(\cX_A,
		\Delta_{\cX_A/\cSTP}
		+\varepsilon\mathcal B_A)
		\longrightarrow\cSTP
		$$
		is a KSBA stable family.
	\end{corollary}
	
	\begin{proof}
		Set
		$D=\Delta_{\cX_A/\cSTP}+\varepsilon\mathcal B_A.$
		The morphism $f$ is flat and projective. Every geometric fiber $(X_s,B_s)$ is a stable toric pair.
      Hence $(X_s,D_s)$ is a KSBA stable pair for sufficiently
      small $\varepsilon>0$ by
      \cite[Lemma~4.4]{Alexeev06}. Indeed, the underlying polyhedral complex is
		a subdivision of the convex polytope $Q$,
		and $B_s$ contains no torus orbit. By Proposition~\ref{prop:relative-log-canonical},
		$
		K_{\cX_A/\cSTP}+D
		\sim_{\QQ}\varepsilon\mathcal B_A,
		$
		which is $\QQ$-Cartier and relatively ample. The boundary divisors are relative Mumford
		divisors. Since $\cSTP$ is reduced, they are
		K-flat by \cite[7.4.2]{Kollar23}.
		Thus the family is stable in the sense of
		\cite[Definition--Theorem~4.7]{Kollar23}.
	\end{proof}

	For this fixed $0<\varepsilon\ll 1$, write
	$
	\Delta_{\cX_A/\cSTP}
	=
	\Delta_1+\cdots+\Delta_s
	$
	for the marking of the relative toric boundary divisors, and set
	$
	\mathbf a_A
	=
	(1,\ldots,1,\varepsilon).
	$
	Thus the associated KSBA pair is regarded as the marked pair
	$$
	\left(
	X,\,
	\Delta_1+\cdots+\Delta_s+\varepsilon B
	\right),
	$$
	where the divisors $\Delta_i$ are marked separately. Following Koll\'ar
	\cite[Definition~8.3 and Theorem~8.1]{Kollar23},
	let
	$
	\mathcal{SP}(\mathbf a_A,d,v_A)
	$
	denote the moduli stack of marked KSBA stable pairs with coefficient
	vector $\mathbf a_A$, dimension $d$, and the corresponding fixed
	volume $v_A$.
	
	By Corollary~\ref{cor:ksba-family}, the universal family determines
	a classifying morphism
	$$
	\Phi_\varepsilon:
	\cSTP
	\longrightarrow
	\mathcal{SP}(\mathbf a_A,d,v_A).
	$$
	
	Let $\mathcal U_A\subset\cSTP$ be the dense open locus parametrizing
	normal stable toric pairs, and define
	$$
	\mathcal{SP}_{A,\mathrm{main}}
	\deq
	\overline{
		\Phi_\varepsilon(\mathcal U_A)
	}_{\mathrm{red}}
	\subset
	\mathcal{SP}(\mathbf a_A,d,v_A).
	$$
	
	\begin{proposition}
		The classifying morphism induces an isomorphism
		$$
		\cSTP
		\longrightarrow
		\bigl(\mathcal{SP}_{A,\mathrm{main}}\bigr)^\nu.
		$$
	\end{proposition}
	
\begin{proof} The classifying morphism is representable, since it induces injections on stabilizers, and proper since $\cSTP$ is proper and the KSBA stack is separated.
	
	To show quasi-finiteness, fix a geometric KSBA pair in the image.
	The marked boundary and conductor determine the full
	toric boundary on each normalization component,
	while $B$ determines its polarization.
	Since $A$ is finite and the component supports are
	full-dimensional, only finitely many compatible
	marked subdivisions and affine lattice identifications
	are possible.
	For each choice, \cite[Theorem~2.6.8]{Alexeev02}
	gives finitely many gluings.
	The geometric fibers are therefore finite, also
	after normalization, so the morphism is finite.
	
	On the normal full-support locus, an isomorphism
	preserving each labeled boundary divisor fixes all
	primitive fan rays, hence has trivial lattice part
	and is $T$-equivariant.
	Thus the morphism is birational.
	Normality of $\cSTP$ then gives a finite birational
	morphism to
	$\bigl(\mathcal{SP}_{A,\mathrm{main}}\bigr)^\nu$,
	which is an isomorphism since the target is normal.
\end{proof}
	
	\subsection{Kappa Classes and Coarse Descent}
	We define the kappa classes which are the main object of this paper. We work with operational Chow rings defined in \cite[Chapter~17]{Ful}, following \cite{FM}. Throughout this paper, we assume that all operational Chow rings have rational coefficients.
	
	Recall the definition of kappa classes on a KSBA moduli stack in \cite{AlexeevKappaBasic}. Fix a KSBA moduli stack $\cM$ parametrizing stable pairs of dimension
	$d$ with fixed numerical invariants, and let
	$
	f:(\cX,\mathcal D)\longrightarrow \cM
	$
	be the universal family. Then $f$ has pure relative dimension $d$, and
	the relative log canonical divisor
	$
	K_{\cX/\cM}+\mathcal D
	$
	defines a relatively ample $\QQ$-line bundle on $\cX$.
	
	Since the universal family $f$ is flat and proper of pure relative
	dimension $d$, the flat orientation in bivariant Chow theory defines
	a generalized Gysin pushforward
	$$
	f_*:A^p(\cX)_\QQ\longrightarrow A^{p-d}(\cM)_\QQ.
	$$
	
	\begin{definition}
		\cite{AlexeevKappaBasic}
		For $r\geq 0$, the $r$-th kappa class is
		$$
		\kp_r
		=
		f_*
		c_1\left(
		\cO_{\cX}(K_{\cX/\cM}+\mathcal D)
		\right)^{d+r}
		\in A^r(\cM)_\QQ.
		$$
	\end{definition}
	
	In our setting for stable toric pairs, we have $\mathcal D = \Delta_{\cX_A/\cSTP}
	+\varepsilon\mathcal B_A$, and 
	$$
	K_{\cX_A/\cSTP} + \cD 
	\sim_{\QQ}
	\varepsilon\mathcal B_A
	$$
	by Proposition~\ref{prop:relative-log-canonical}.
	We write $\kp_r(\varepsilon)$ for the kappa class of this toric pair,
	making its dependence on $\varepsilon$ explicit. Hence
	$$
	\kp_r(\varepsilon)
	=
	\varepsilon^{d+r}
	{f}_*c_1(\mathcal B_A)^{d+r} \in A^r(\cSTP)_\QQ.
	$$
	This motivates the following normalization.
	\begin{definition}
		The $r$-th normalized kappa class is 
		$$
		\kp_r
		\deq
		{f}_*c_1(\mathcal B_A)^{d+r} \in A^r(\cSTP)_\QQ.
		$$
	\end{definition}
	One can define kappa classes and their normalized versions on the coarse moduli space. Let $
	\rho:
	\cSTP
	\longrightarrow
	\STP
	$
	be the coarse morphism. By
	\cite[Proposition~6.1]{Vis}, pullback along $\rho$ induces an
	isomorphism
	$$
	\rho^*:
	A^*(\STP)_\QQ
	\xrightarrow{\sim}
	A^*(\cSTP)_\QQ.
	$$
	Therefore, by \cite{AlexeevKappaBasic}, there exist unique descended classes
	$$
	\bkp_r(\varepsilon),\ \bkp_r
	\in
	A^r(\STP)_\QQ
	$$
	such that
	$
	\rho^*\bkp_r(\varepsilon)=\kp_r(\varepsilon),
	$ and $\rho^*\bkp_r=\kp_r.
	$
	We call $\bkp_r(\varepsilon)$ and $\bkp_r$ the coarse descendants of
	$\kp_r(\varepsilon)$ and $\kp_r$, respectively.

	\section{Toric Models for Stable Toric Pairs}\label{sec:toric-models}
	In this section, we construct an explicit toric family over the moduli space $\STP$ and modify lattices to build a flat family with reduced fibers used to compute the kappa classes.
	
	\subsection{The Effective Toric Model}
	
	We describe an explicit toric family over the moduli space $\STP$ using the fiber-polytope construction in \cite[Section~4.7]{AlexeevPardini}. For a more detailed account of fiber polytopes, see \cite{BilleraSturmfels92}. Related toric families over spaces parametrizing matroidal
	subdivisions were constructed by Lafforgue
	\cite{Lafforgue2003}. 
	
	Let
	$
	\psi:M_P\to M_Q
	$
	be a surjective map of lattices and let $P\subset M_{P,\RR}$ be a
	lattice polytope with image $Q=\psi(P)$. The associated fiber polytope
	$\Sigma(P,Q)$ determines an equidimensional toric morphism
	$$
	X_{P+\Sigma(P,Q)}
	\longrightarrow
	X_{\Sigma(P,Q)}
	$$
	where $+$ denotes Minkowski addition. The strata of the base correspond to coherent subdivisions of $Q$ induced by $\psi$, and
	the reduced irreducible components of the corresponding fibers are
	indexed by the maximal cells of these subdivisions.
	
	We apply this construction to the fixed $A \subset M$ and $Q = \operatorname{Conv}(A) \subset M_{\RR}$. Let $P$ be the standard simplex $\Delta_A=\operatorname{Conv}\{e_a\mid a\in A\}$ and $\PP^A$ be the corresponding projective space. The translation lattice of $\Delta_A$ is
	$$
	M_{\PP^A}
	\deq
	\{
	(u_a)_{a\in A}\in\ZZ^A
	|
	\sum_{a\in A}u_a=0
	\}.
	$$
	The configuration $A$ defines a homomorphism
	$$
	\phi:M_{\PP^A}\longrightarrow M,
	\qquad
	(u_a)_{a\in A}\longmapsto\sum_{a\in A}u_a a.
	$$
	Set
	$\Lambda_A\deq\operatorname{im}(\phi)=\ZZ(A-A)$
	and
	$$
	L_A\deq\ker(\phi)
	=
	\{
	(u_a)_{a\in A}\in\ZZ^A
	|
	\sum_{a\in A}u_a=0,\quad
	\sum_{a\in A}u_a a=0
	\}.
	$$
	Since $A$ is full-dimensional, $\Lambda_A$ has finite index in $M$. Thus
	\begin{equation}\label{exact}
		0
		\longrightarrow
		L_A
		\longrightarrow
		M_{\PP^A}
		\longrightarrow
		\Lambda_A
		\longrightarrow
		0
	\end{equation}
	is exact. We call $\Lambda_A$ the
	\emph{effective character lattice associated with $A$}.
	The map $\phi$ is surjective onto $\Lambda_A$, and we set
	$g_A\deq[M:\Lambda_A]$. Since
	$\Lambda_{A,\RR}=M_{\RR},$
	the induced map of real vector spaces remains the same. Thus the affine map
	$$
	\Delta_A\longrightarrow Q=\operatorname{Conv}(A),
	\qquad
	e_a\longmapsto a,
	$$
	gives the fiber polytope $\Sigma(\Delta_A,Q)$, which has the same normal fan as the secondary polytope $\Sec(A)$. The Minkowski sum construction yields a toric morphism
	$$
	\pi^{\mathrm{eff}}:
	X^{\mathrm{eff}}
	\longrightarrow
	\STP,
	$$
	where \mbox{$X^{\mathrm{eff}}\deq X_{\Delta_A+\Sec(A)}$}.
	Since the normal fan of $\Delta_A+\Sec(A)$ refines the normal fan of $\Delta_A$,
	there is also a birational toric morphism
	$$
	p^{\mathrm{eff}}:
	X^{\mathrm{eff}}
	\longrightarrow
	X_{\Delta_A}=\PP^A.
	$$
Thus $\pi^{\mathrm{eff}}$ is an equidimensional toric morphism
with generic fiber character lattice~$\Lambda_A$.
The reduced components of each fiber are indexed by the
maximal cells of the corresponding subdivision.
	$$
	\begin{tikzcd}[column sep=3em]
		& X^{\mathrm{eff}}
		\arrow[dl,"p^{\mathrm{eff}}"']
		\arrow[dr,"\pi^{\mathrm{eff}}"] & \\
		\PP^A && \STP 
	\end{tikzcd}
	$$

	The effective family need not be flat, and its
	fibers need not be reduced. Moreover, its
	generic fiber has character lattice
	$\Lambda_A$, which may differ from $M$.
	
	The following example describes a nonreduced
	fiber both by an equation and by the
	corresponding map of normal fans.
	
\begin{example}\label{ex:three-points}
Let $A=\{0,1,2\}\subset\mathbb Z$.
There are two regular triangulations of $A$:
$$
\mathcal T_{\mathrm{split}}
=
\bigl\{\{0,1\},\{1,2\}\bigr\},
\qquad
\mathcal T_{\mathrm{omit}}
=
\bigl\{\{0,2\}\bigr\}.
$$
In the second triangulation, the interior point $1$
is omitted. Their GKZ vectors are
$$
\varphi_{\mathcal T_{\mathrm{split}}}
=(1,2,1),
\qquad
\varphi_{\mathcal T_{\mathrm{omit}}}
=(2,0,2).
$$

Thus $\Sec(A)$ is a line segment and the coarse moduli space is
$\STP\simeq\mathbb P^1.$ The two torus-fixed points correspond to the two triangulations. For $\mathcal T_{\mathrm{omit}}$, the underlying curve is $\PP^1$, and its polarization has degree~$2$.
For $\mathcal T_{\mathrm{split}}$, the underlying curve
consists of two copies of $\PP^1$ meeting at one node.

Now we describe the effective model of this example. On a suitable affine chart of the base
	$\STP\simeq\mathbb P^1$, the effective toric
	family is given by
	$$
	z_1^2=t z_0z_2
	$$
	inside $\mathbb P^2\times\mathbb A^1_t$. For $t\neq0$, the fiber is a smooth conic,
whereas the fiber over $t=0$ is a nonreduced double line.

Figure~\ref{fig:three-points-effective-fan}
describes the effective model in terms of fans.
The fan $\Sigma_{X^{\mathrm{eff}}}$ refines~$\Sigma_{\PP^2}$, and the identity on their common
lattice induces $p^{\mathrm{eff}}$, as indicated
by the left arrow. The right arrow is the quotient map $\beta$
defining $\pi^{\mathrm{eff}}$.
In the upper panel, this map is represented by
orthogonal projection onto the red line,
which carries the secondary fan.
The primitive generator of the indicated source ray
maps to $2u$, where $u$ is the primitive generator
of the corresponding base ray.
Thus the associated toric divisor occurs with
multiplicity $2$ in the fiber over $t=0$,
recovering the double fiber directly from the fan map.

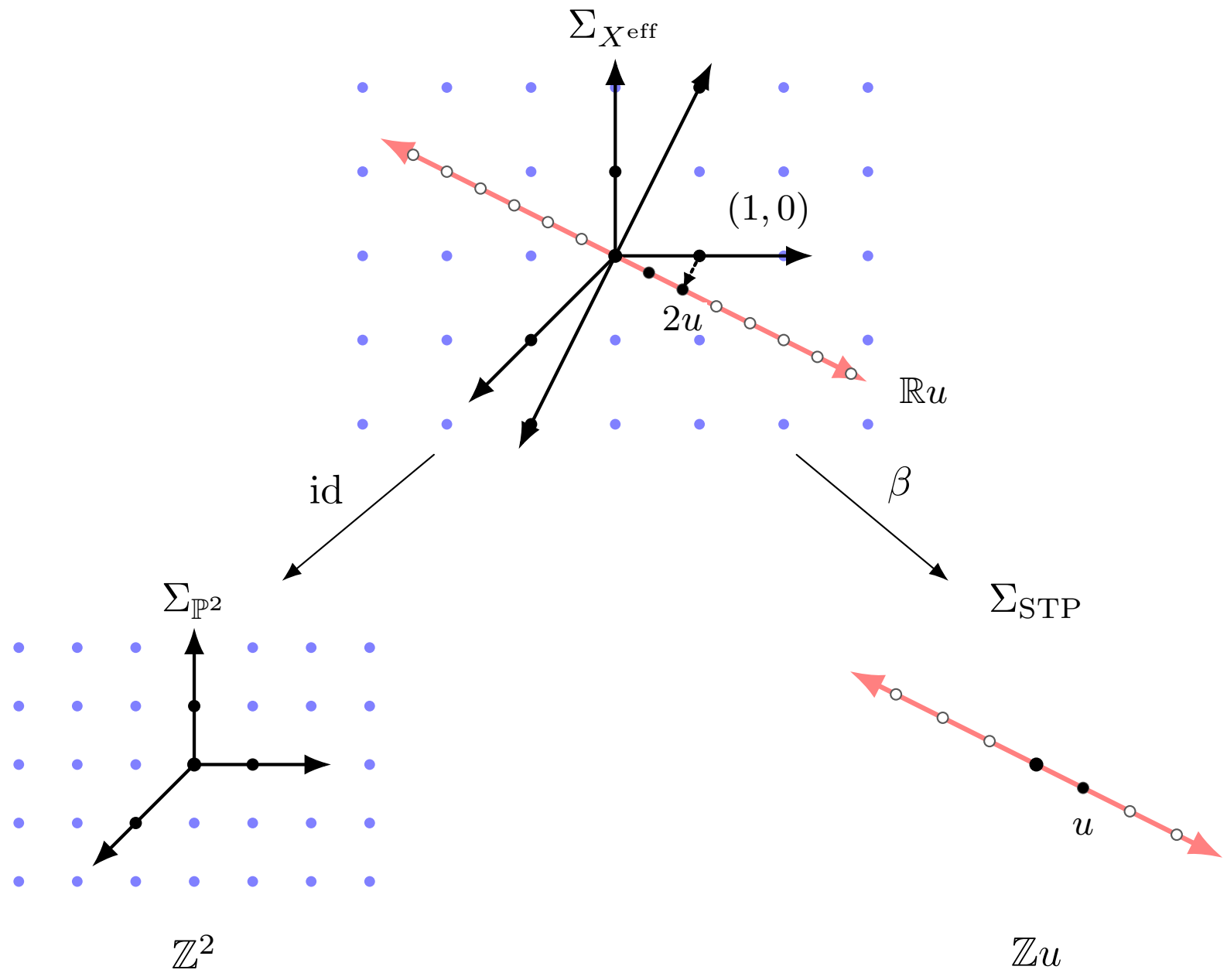
\begin{figure}[htbp]
	\centering
	\begin{tikzpicture}[
		x=1cm,y=1cm,
		>=Latex,
		line cap=round,
		line join=round
		]
		
		\node at (0,3.95)
		{$\Sigma_{X^{\mathrm{eff}}}$};
		
		\begin{scope}[
			shift={(0,2.0)},
			x=0.72cm,y=0.72cm
			]
			\foreach \i in {-3,-2,...,3}{
				\foreach \j in {-2,-1,...,2}{
					\fill[blue!50]
					(\i,\j) circle (1.3pt);
				}
			}
			
			\coordinate (O) at (0,0);
			
			\draw[thick,->] (O) -- (2.35,0);
			\draw[thick,->] (O) -- (0,2.35);
			\draw[thick,->] (O) -- (1.15,2.30);
			\draw[thick,->] (O) -- (-1.75,-1.75);
			\draw[thick,->] (O) -- (-1.15,-2.30);
			
			\foreach \P in {
				(1,0),(0,1),(1,2),(-1,-1),(-1,-2)
			}{
				\fill \P circle (1.5pt);
			}
			
			\draw[red!50,very thick,<->]
			(-2.8,1.4) -- (3.0,-1.5);
			
			\foreach \k in {-6,-5,...,7}{
				\fill[white]
				({0.4*\k},{-0.2*\k})
				circle (1.3pt);
				\draw[black!65]
				({0.4*\k},{-0.2*\k})
				circle (1.3pt);
			}
			
			\coordinate (u) at (0.4,-0.2);
			\coordinate (twou) at (0.8,-0.4);
			\coordinate (source) at (1,0);
			
			\fill (O) circle (1.7pt);
			\fill (u) circle (1.3pt);
			\fill (twou) circle (1.3pt);
			\fill (source) circle (1.5pt);
			
			\draw[
			densely dotted,
			thick,
			-{Latex[length=1.4mm,width=1.2mm]}
			] (source) -- (twou);
			
			\node[above right=3pt,font=\small]
			at (source) {$(1,0)$};
			
			\node[
			below=3pt,font=\small,
			fill=white,inner sep=1pt
			] at (twou) {$2u$};
			
			\node[anchor=west,font=\small]
			at (3.2,-1.6) {$\mathbb R u$};
		\end{scope}
		
		\node (Pfan) at (-3.6,-0.95)
		{$\Sigma_{\mathbb P^2}$};
		
		\node (Bfan) at (3.6,-0.95)
		{$\Sigma_{\STP}$};
		
		\draw[->]
		(-1.55,0.30) --
		node[above left] {$\mathrm{id}$}
		(-2.85,-0.78);
		
		\draw[->]
		(1.55,0.30) --
		node[above right] {$\beta$}
		(2.85,-0.78);
		
		\begin{scope}[
			shift={(-3.6,-2.35)},
			x=0.50cm,y=0.50cm
			]
			\foreach \i in {-3,-2,...,3}{
				\foreach \j in {-2,-1,...,2}{
					\fill[blue!50]
					(\i,\j) circle (1.3pt);
				}
			}
			
			\draw[thick,->] (0,0) -- (2.35,0);
			\draw[thick,->] (0,0) -- (0,2.35);
			\draw[thick,->] (0,0) -- (-1.75,-1.75);
			
			\foreach \P in {
				(1,0),(0,1),(-1,-1)
			}{
				\fill \P circle (1.5pt);
			}
			
			\fill (0,0) circle (1.7pt);
		\end{scope}
		
		\begin{scope}[
			shift={(3.6,-2.35)}
			]
			\draw[red!50,very thick,<->]
			(-1.6,0.8) -- (1.6,-0.8);
			
			\foreach \k in {-3,-2,...,3}{
				\fill[white]
				({0.4*\k},{-0.2*\k})
				circle (1.3pt);
				\draw[black!65]
				({0.4*\k},{-0.2*\k})
				circle (1.3pt);
			}
			
			\fill (0,0) circle (1.7pt);
			\fill (0.4,-0.2) circle (1.3pt);
			
			\node[below=4pt,font=\small]
			at (0.4,-0.2) {$u$};
		\end{scope}
		
		\node at (-3.6,-3.95) {$\mathbb Z^2$};
		\node at (3.6,-3.95) {$\mathbb Z u$};
		
	\end{tikzpicture}
	
	\caption{
		Fan maps of the effective model for
		$A=\{0,1,2\}$.
		The image $2u$ of the primitive generator of a source ray
		records the multiplicity of the nonreduced
		fiber. The blue dots represent points of
		the ambient lattice.
	}
	\label{fig:three-points-effective-fan}
\end{figure}

Replacing the base lattice $\ZZ u$ by $2\ZZ u$
corresponds to the base change $t=s^2$.
The normalized pullback has source lattice
$\beta^{-1}(2\ZZ u)$.
The real fans remain unchanged, and the indicated
primitive source vector now maps to the primitive
base vector~$2u$, yielding a reduced fiber.

In this example, the effective family is flat,
but its special fiber is not reduced. More
explicitly, after the base change $t=s^2$,
the normalization is described by
$$
([u:v],s)\longmapsto
([u^2:suv:v^2],s).
$$
The central fiber is the reduced curve
$\mathbb P^1$ corresponding to the triangulation
$\mathcal T_{\mathrm{omit}}$ introduced above. Its map to the
reduced line $\{z_1=0\}$ has degree~$2$.
\end{example}

In the next subsection, we generalize this construction
to arbitrary configurations, modifying the lattices of the
effective model to obtain a flat family with reduced fibers.
	
	\subsection{Kummer Model and Lattice Modification}\label{subsec:kummer-model}
	We now modify the lattices of the effective toric model to obtain a flat family with reduced fibers and generic fiber character lattice $M$. We use the weak semistability criterion of \cite{Molcho}.
	
	Dualizing \eqref{exact}
	gives an exact sequence
	$$
	0
	\longrightarrow
	N_{\Lambda_A}
	\longrightarrow
	N_{\PP^A}
	\overset{\beta}{\longrightarrow}
	N_{\Sec(A)}
	\longrightarrow
	0.
	$$
	Here $N$, $N_{\Lambda_A}$, $N_{\PP^A}$, and $N_{\Sec(A)}$ denote the dual lattices of $M$, $\Lambda_A$, $M_{\PP^A}$, and $L_A$, respectively.
	The inclusion $\Lambda_A\subset M$ induces an inclusion
	$
	N \subset N_{\Lambda_A},
	$
	of finite index
	$[N_{\Lambda_A}:N]=g_A.$
	In particular, $g_A N_{\Lambda_A}\subseteq N.$ 
	
	Note that, under the pairing
	$$
	\langle[\omega],u\rangle
	=
	\sum_{a\in A}\omega(a)u_a,
	$$
	we identify $N_{\Sec(A),\RR}$ with the space of real height functions
	on $A$ modulo affine functions.
	If $\tau_{\cS}$ is the cone corresponding to a regular
	subdivision $\cS$, representatives of its classes are
	affine on every cell of $\cS$.
	
	We first choose a common multiple that removes all lattice defects in
	the cone maps of the effective family. Since $\beta$ defines an equidimensional toric morphism, any cone $\sigma\in\Sigma_{X^{\mathrm{eff}}}$ maps onto a cone,  say $\tau=\beta(\sigma)\in\Sigma_{\Sec(A)}$.
	For a cone $\sigma\in\Sigma_{X^{\mathrm{eff}}}$ and its image
	$\tau=\beta(\sigma)$, set
	$
	N_\sigma
	\deq
	N_{\PP^A}\cap\operatorname{Span}_{\RR}(\sigma),
	$ and $
	N_\tau
	\deq
	N_{\Sec(A)}\cap\operatorname{Span}_{\RR}(\tau).
	$
	For each such pair $(\sigma,\tau)$ with $\dim\sigma=\dim\tau$, define
	$$
	m_{\sigma/\tau}
	\deq
	[N_\tau:\beta(N_\sigma)]
	$$
	and 
	$$
	m
	\deq
	\operatorname{lcm}
	\left\{
	m_{\sigma/\tau}
	\ \middle|\
	\beta(\sigma)=\tau,\ 
	\dim\sigma=\dim\tau
	\right\}.
	$$
	By construction,
	$mN_\tau
	\subseteq
	\beta(N_\sigma)$
	for every such pair $(\sigma,\tau)$. 
	
	The factor $m$ will ensure reduced fibers, but an additional factor is needed to recover the original character lattice $M$. We therefore set~\mbox{$\mathfrak m\deq g_A m$}.
	Define the sublattices
	$$
	\widetilde N_{\Sec(A)}
	\deq
	\mathfrak m N_{\Sec(A)} \subset N_{\Sec(A)}
	$$
	and
	$$
	\widetilde N_X
	\deq
	N+\mathfrak m N_{\PP^A}
	\subset N_{\PP^A}.
	$$ 
	The restriction of $\beta$ to $\widetilde N_X$ gives an exact
	sequence
	$$
	0
	\longrightarrow
	N
	\longrightarrow
	\widetilde N_X
	\overset{\widetilde\beta}{\longrightarrow}
	\widetilde N_{\Sec(A)}
	\longrightarrow
	0.
	$$
	Indeed,
	$$
	\beta(\widetilde N_X)
	=
	\beta(N+\mathfrak m N_{\PP^A})
	=
	\mathfrak m N_{\Sec(A)},
	$$
	while
	$$
	\begin{aligned}
		\ker(\widetilde\beta)
		&=
		\widetilde N_X\cap N_{\Lambda_A}\\
		&=
		N+\mathfrak m N_{\Lambda_A}\\
		&=
		N.
	\end{aligned}
	$$
	The last equality follows from
	$g_A N_{\Lambda_A}\subseteq N$
	and $g_A\mid\mathfrak m$.
	
	Keeping the same real fans, define
	$$
	\tSTP
	\deq
	X\bigl(
	\Sigma_{\Sec(A)},
	\widetilde N_{\Sec(A)}
	\bigr)
	$$
	and
	$$
	\widetilde X
	\deq
	X\bigl(
	\Sigma_{X^{\mathrm{eff}}},
	\widetilde N_X
	\bigr).
	$$
	The lattice homomorphism
	$
	\widetilde\beta:
	\widetilde N_X
	\longrightarrow
	\widetilde N_{\Sec(A)}
	$ induces a toric morphism
	\mbox{$\widetilde\pi:\widetilde X\longrightarrow\tSTP$}. Similarly, the inclusion
	$
	\widetilde N_X\subset N_{\PP^A},
	$ together with the fan refinement
	$\Sigma_{X^{\mathrm{eff}}}\to\Sigma_{\PP^A}$,
	induces a toric morphism
	\mbox{$\widetilde p:\widetilde X\longrightarrow\PP^A$}.
	We also have a finite toric morphism
	$
	\mu:
	\tSTP\rightarrow \STP
	$,
	and the diagram
	$$
	\begin{tikzcd}[column sep=4em, row sep=2.5em]
		\widetilde X
		\arrow[r]
		\arrow[d, "\widetilde\pi"']
		\arrow[dr, phantom, "\circlearrowright" description]
		&
		X^{\mathrm{eff}}
		\arrow[d, "\pi^{\mathrm{eff}}"]
		\\
		\tSTP
		\arrow[r, "\mu"']
		&
		\STP
	\end{tikzcd}
	$$
	commutes, but in general it is not Cartesian.

	When $g_A=1$, this construction is the
	normalized pullback of the effective model
	along $\mu$. For general $g_A$, an additional
	finite lattice modification changes the
	character lattice of the generic fiber from
	$\Lambda_A$ to $M$. Thus the construction has
	two purposes: removing the lattice indices
	responsible for nonreduced fibers and
	recovering the prescribed fiber lattice.

	\begin{lemma}\label{lem:face-isomorphic-to-target}
		Let $\sigma\in\Sigma_{X^{\mathrm{eff}}}$ with $\beta(\sigma)=\tau$. Then there is a
		face $\sigma'\preceq\sigma$ such that
		$
		\beta_{\mathbb R}\colon \sigma'\xrightarrow{\sim}\tau.
		$	
	\end{lemma}
	
	\begin{proof}
		
		Choose a point $w\in\operatorname{relint}(\tau)$
		and a vertex $v$ of
		$\sigma\cap\beta_{\mathbb R}^{-1}(w)$.
		Let $\sigma'$ be the smallest face of $\sigma$
		containing $v$. By equidimensionality, $\beta_{\mathbb R}(\sigma')$
		is a cone contained in
		$\beta_{\mathbb R}(\sigma)=\tau$.
		Since $w=\beta_{\mathbb R}(v)$ belongs to
		$\beta_{\mathbb R}(\sigma')$ and lies in
		$\operatorname{relint}(\tau)$, we conclude that
		$
		\beta_{\mathbb R}(\sigma')=\tau.
		$
		The injectivity of the restriction of $\beta_\RR$ to $\Span_\RR(\sigma')$ follows from the fact that $v$ is a vertex.
	\end{proof}
	
	\begin{proposition}\label{prop:weakly-semistable}
		The toric morphism
		$$
		\widetilde\pi:
		\widetilde X\longrightarrow\tSTP
		$$
		is weakly semistable.
	\end{proposition}
	
	\begin{proof}
		By \cite[Definition~2.1.2]{Molcho}, it suffices to show that every cone maps onto a cone and that the induced map on lattice points is surjective. 
		
		Since $\pi^{\mathrm{eff}}$ is equidimensional, every cone
		$\sigma\in\Sigma_{X^{\mathrm{eff}}}$ maps surjectively onto a cone
		$\tau\in\Sigma_{\Sec(A)}$. It remains to prove
		$$
		\widetilde\beta(\widetilde N_X\cap\sigma)
		=
		\widetilde N_{\Sec(A)}\cap\tau.
		$$
		
		The inclusion from left to right is immediate. Conversely, let
		$$
		y\in
		\widetilde N_{\Sec(A)}\cap\tau
		=
		\mathfrak mN_{\Sec(A)}\cap\tau.
		$$
		Write
		$y=\mathfrak mz=g_Amz
		$ for some $z\in N_\tau\cap\tau.
		$
		By Lemma~\ref{lem:face-isomorphic-to-target}, there is a face
		$\sigma'\preceq\sigma$ such that
		$\beta_{\RR}\colon\sigma'\xrightarrow{\sim}\tau.$ Since $m_{\sigma'/\tau}\mid m$, we have $mN_\tau\subseteq\beta(N_{\sigma'}).$
		Thus there is a unique
		$
		x'\in N_{\sigma'}\cap\sigma'
		$
		such that
		$
		\beta(x')=mz.
		$	
		
		We claim that $x \deq g_Ax' \in \widetilde{N}_X \cap \sigma$ and $\widetilde \beta(x) = y.$
		Since $\beta(x')=mz$, we have
		$$
		x'\in
		\beta^{-1}(mN_{\Sec(A)})
		=
		N_{\Lambda_A}+mN_{\PP^A}.
		$$
		Therefore
		$$
		g_Ax'
		\in
		g_AN_{\Lambda_A}+g_AmN_{\PP^A}
		\subseteq
		N+\mathfrak mN_{\PP^A}
		=
		\widetilde N_X.
		$$
		Since $x'\in\sigma$, we also have $g_Ax'\in\sigma$.  Finally, $\widetilde\beta(x)=g_Amz=y$.
	\end{proof}
	
	\begin{corollary}\label{cor:flat-reduced}
		The morphism
		$$
		\widetilde\pi:
		\widetilde X
		\longrightarrow
		\tSTP
		$$
		is flat with reduced fibers. The fibers are stable toric
		varieties with character lattice $M$ and subdivision data
		determined by $A$.
	\end{corollary}
	
	\begin{proof}
		Flatness and reducedness of the fibers follow from  \cite[Theorem~2.1.4]{Molcho}, which states that a toric morphism
		is weakly semistable if and only if it is flat with reduced fibers.
		
		Changing the lattice does not change the underlying real fans,
		so the fibers have the same subdivision data as those of
		$\pi^{\mathrm{eff}}$. Moreover,
		$\ker(\widetilde\beta)=N$, so the fiber torus has character
		lattice $N^\vee = M$.
	\end{proof}

	We call these finite lattice changes the
	\emph{Kummer modification} and the resulting family $\widetilde{\pi}$ the {\it Kummer model}. For a maximal cell $C$, the corresponding fiber component
	has polarization polytope $\operatorname{Conv}(C)$ in the
	lattice~$M$. The restriction of $\widetilde p$ to each fiber
	is finite, so $\widetilde p^*\cO_{\PP^A}(1)$ is relatively
	ample. Moreover, the pullback of
	$$
	H_A=\left\{\sum_{a\in A}z_a=0\right\}\subset\PP^A
	$$
	contains no torus orbit. Thus $\widetilde{\mathcal B}:=\widetilde p^*H_A$
	defines a family of stable toric pairs associated to $A$
	over $\tSTP$.
		
	\subsection{Kummer Kappa Classes}
	The Kummer family $\widetilde{\pi}$ determines a classifying morphism
	\mbox{$\chi:\tSTP\longrightarrow\cSTP$}.
	By construction,
	$
	\mathcal O_{\widetilde X}(\widetilde{\mathcal B})
	\cong
	\widetilde p^*\mathcal O_{\PP^A}(1).
	$
By the base-change compatibility of kappa classes
\cite{AlexeevKappaBasic}, their pullbacks to the Kummer
base can be computed directly from this family.
	
	\begin{definition} The $r$-th Kummer kappa class is defined by the pullback
		$$
		\tbkp_r
		\deq
		\chi^*\kp_r
		=
		\widetilde\pi_*
		\left(
		(\widetilde p^*h)^{d+r}
		\right)
		\in
		A^r(\tSTP)_\QQ
		$$
		where $h = c_1(\mathcal O_{\PP^A}(1)).$
	\end{definition}

The following proposition identifies the Kummer kappa classes
with pullbacks of the coarse classes.
	\begin{proposition}\label{prop:kummer-kappa-pullback}
		The pulled-back kappa class satisfies
		$$
		\tbkp_r
		=
		\mu^*\bkp_r.
		$$
	\end{proposition}
	\begin{proof}
		By the moduli interpretation, for every geometric
		point $s\in\tSTP$, the points $\mu(s)$ and
		$\rho(\chi(s))$ represent the same isomorphism
		class of stable toric pairs. Thus the two morphisms
		agree on geometric points.
		
		Consider their graphs
		\mbox{$\Gamma_\mu,\ \Gamma_{\rho\circ\chi}\subset\tSTP\times\STP$}.
		Since $\STP$ is separated, both graphs are
		closed subschemes. Moreover, each graph is
		isomorphic to $\tSTP$, which is normal and
		hence reduced. The two graphs have the same geometric points,
		so they coincide as reduced closed subschemes.
		Consequently,
		$
		\Gamma_\mu=\Gamma_{\rho\circ\chi},
		$
		and hence~\mbox{$\mu=\rho\circ\chi$}.
		
		By functoriality of pullback, 
		$$
		\tbkp_r
		=\chi^*\kp_r
		=\chi^*\rho^*\bkp_r
		=\mu^*\bkp_r.
		$$
	\end{proof}
	
	\section{The Main Intersection Formula}\label{sec:intrinsic-kummer}
	
	In this section, we determine the normalized kappa classes using
	the Kummer family constructed in Section~\ref{sec:toric-models}.
	The computation relies on the combinatorial geometry of the
	torus-invariant strata of $\STP$.
	
	Since $\STP$ is a complete toric variety, by
	Fulton--Sturmfels, we obtain an isomorphism
	$$
	A^r(\STP)_\QQ
	\xrightarrow{\sim}
	\Hom\bigl(A_r(\STP)_\QQ,\QQ\bigr),
	\qquad
	\alpha\longmapsto
	\bigl([Z]\mapsto\deg(\alpha\cap[Z])\bigr).
	$$
	Moreover, $A_r(\STP)$ is generated by the $r$-dimensional
	torus-invariant subvarieties. Thus $\bkp_r$ is determined
	by its degrees on the $r$-dimensional orbit closures
	$V_{\cS}\subset\STP$.

	\subsection{Invariant Orbit Closures and Lattice Indices}\label{subsec:common-invariant-strata}
    We first work with the effective model, since the Kummer modification preserves the real fans. We first identify the contributing components and compute
	the lattice indices that appear in the intersection formula.
	The results in this section will also be used in the direct toric computation of Section~\ref{sec:alternative-toric}.
	
Fix a regular subdivision $\cS$ of $(Q,A)$, and let
$\tau_{\cS}\in\Sigma_{\Sec(A)}$ correspond to the
$r$-dimensional orbit closure $V_{\cS}\subset\STP$.
Set
$X_{\cS}\deq
(\pi^{\mathrm{eff}})^{-1}(V_{\cS}).$ For each maximal cell $C\in\cS^{\max}$, let
$\sigma_{\cS,C}\in\Sigma_{X^{\mathrm{eff}}}$ be the cone
corresponding to the component indexed by $C$ over
$O(\tau_{\cS})$.
Set
$X_{\cS,C}:=V(\sigma_{\cS,C})\subset X_{\cS}$.
As $C$ ranges over $\cS^{\max}$, these are the
irreducible components of $X_{\cS}$, each endowed
with its reduced induced structure.

Let $\gamma_C\in\Sigma_{\PP^A}$ be the cone corresponding
to the coordinate subspace $\PP^C\subset\PP^A$.
The restriction of the polarization morphism is denoted by
$$
p_{\cS,C}:X_{\cS,C}\longrightarrow
\PP^C\hookrightarrow\PP^A.
$$

We apply the lattice construction for $A$ to each
maximal cell $C$: set
$$
\Lambda_C\deq\ZZ(C-C),\qquad
g_C\deq[M:\Lambda_C],\qquad
M_{\PP^C}\deq
\left\{(u_a)\in\ZZ^C\ \middle|\ \sum_{a\in C}u_a=0\right\},
$$
and define $L_C$ by the exact sequence
\begin{equation}\label{exactc}
	0\longrightarrow L_C\longrightarrow M_{\PP^C}
	\xrightarrow{\phi_C}\Lambda_C\longrightarrow0,
	\qquad
	\phi_C((u_a))=\sum_{a\in C}u_a a.
\end{equation}
	We first determine which maximal cells can contribute to an intersection
	with the $r$-th kappa class.
	
	\begin{lemma}\label{lem:cell-size-bound}
		Let $V_{\cS}\subset\STP$ be an $r$-dimensional torus-invariant
		subvariety and $C\in\cS^{\max}$. Then
		$|C|\leq d+r+1.$
	\end{lemma}

	\begin{proof}
		One can regard $L_C$ as a
		sublattice of $L_A$ using extension by $0$ outside $C$. Every height function in $\tau_{\cS}$ is
		affine on $C$. Hence $L_{C,\RR}$ annihilates
		$\Span_{\RR}(\tau_{\cS})$. Therefore
		$
		|C|-d-1
		=
		\rk L_C
		\leq
		\rk L_A-\dim\tau_{\cS}
		=r.$
		Thus $|C|\leq d+r+1$.
	\end{proof}

	\begin{corollary}\label{cor:vanishing-nonmaximal-cells}
		If $|C|<d+r+1$, then $p_{\cS,C}^*h^{d+r}=0.$
	\end{corollary}
	
	\begin{proof}
		Since $p_{\cS,C}$ factors through $\PP^C$ and
		$\dim\PP^C=|C|-1<d+r,$
		we have $h^{d+r}|_{\PP^C}=0$. Therefore
		$p_{\cS,C}^*h^{d+r}=0$.
	\end{proof}
	
	Consequently, only the maximal cells $C\in\cS^{\max}$ satisfying
	$|C|=d+r+1$
	can contribute to the intersection numbers below.
	
	\begin{lemma}\label{lem:sigma-coordinate-cone}
		Let $C\in\cS^{\max}$ with $|C|=d+r+1$. Then
		$
		N_{\sigma_{\cS,C}}=N_{\gamma_C}.$
		Consequently,
		$
		N_{\PP^C}
		=
		N_{\PP^A}/N_{\sigma_{\cS,C}}.
		$
	\end{lemma}
	
	\begin{proof}
		The linear subspace $\PP^C\subset\PP^A$ corresponds to the cone
		generated by the rays indexed by $A\setminus C$. The cone
		$\sigma_{\cS,C}$ is contained in this cone. Moreover,
		$$
		\dim\sigma_{\cS,C}
		=
		\dim\tau_{\cS}
		=
		|A|-|C|,
		$$
		so the two cones have the same linear span. Hence their associated
		sublattices of~$N_{\PP^A}$ coincide, and the quotient description of
		$N_{\PP^C}$ follows.
	\end{proof}
	
	\begin{lemma}\label{lem:local-effective-index}
		Let $C\in\cS^{\max}$ satisfy $|C|=d+r+1$. Then
		$$
		m_{\sigma_{\cS,C}/\tau_{\cS}}
		=
		[\Lambda_A:\Lambda_C].
		$$
	\end{lemma}
	
	\begin{proof}
By Lemma~\ref{lem:sigma-coordinate-cone}, the natural
projection gives an exact sequence
$$
0\longrightarrow N_{\sigma_{\cS,C}}
\longrightarrow N_{\PP^A}
\longrightarrow N_{\PP^C}
\longrightarrow 0.
$$
Since $L_C=L_A\cap M_{\PP^C}$ is saturated in $L_A$,
dualizing the inclusion gives a surjective restriction map
$$
\operatorname{res}_C:
N_{\Sec(A)}\longrightarrow N_{\Sec(C)}.
$$

Height functions in $\tau_{\cS}$ are affine on $C$
and therefore annihilate $L_C$. Thus
$$
\Span_{\RR}(\tau_{\cS})
\subseteq \ker(\operatorname{res}_{C, \RR}).
$$
Since $\rk N_{\Sec(C)}=|C|-d-1=r$,
surjectivity gives
$$
\dim\ker(\operatorname{res}_{C, \RR})
=
\rk N_{\Sec(A)}-r
=
\dim\tau_{\cS}.
$$
The inclusion is therefore an equality.
Intersecting with $N_{\Sec(A)}$, we conclude that~\mbox{$N_{\tau_{\cS}}=\ker(\operatorname{res}_C)$}.
		
		Consider the following commutative diagram:
		$$
		\begin{tikzcd}[column sep=3em,row sep=2.3em]
			& 0 \arrow[r] \arrow[d]
			& N_{\sigma_{\cS,C}} \arrow[r,"\beta"] \arrow[d,hook]
			& N_{\tau_{\cS}} \arrow[d,hook]
			& {} \\
			0 \arrow[r]
			& N_{\Lambda_A} \arrow[r] \arrow[d,hook]
			& N_{\PP^A} \arrow[r,"\beta"] \arrow[d,two heads]
			& N_{\Sec(A)} \arrow[r] \arrow[d,two heads]
			& 0 \\
			0 \arrow[r]
			& N_{\Lambda_C} \arrow[r] \arrow[d,two heads]
			& N_{\PP^C} \arrow[r] \arrow[d]
			& N_{\Sec(C)} \arrow[r] \arrow[d]
			& 0 \\
			& N_{\Lambda_C}/N_{\Lambda_A} \arrow[d]
			& 0 \arrow[d]
			& 0 \arrow[d]
			& {} \\
			& 0
			& 0
			& 0
			& {}
		\end{tikzcd}
		$$
		
		The Snake Lemma gives an exact sequence
		$$
		0
		\longrightarrow
		N_{\sigma_{\cS,C}}
		\overset{\beta}{\longrightarrow}
		N_{\tau_{\cS}}
		\longrightarrow
		N_{\Lambda_C}/N_{\Lambda_A}
		\longrightarrow
		0.
		$$
		
		Therefore
		$$
		\begin{aligned}
			m_{\sigma_{\cS,C}/\tau_{\cS}}
			&=
			[N_{\tau_{\cS}}:
			\beta(N_{\sigma_{\cS,C}})]
			\\
			&=
			|N_{\Lambda_C}/N_{\Lambda_A}|
			\\
			&=
			[\Lambda_A:\Lambda_C].
		\end{aligned}
		$$
	\end{proof}
    
	To summarize, only maximal cells satisfying $|C|=d+r+1$ can contribute to the intersection numbers,
	and for each such cell the local lattice index is
	$$
	m_{\sigma_{\mathcal S,C}/\tau_{\mathcal S}}
	=
	[\Lambda_A:\Lambda_C]
	=
	\frac{g_C}{g_A}.
	$$

	\subsection{Computation on the Kummer Family}
	Recall the Kummer family
	\mbox{$\widetilde\pi:\widetilde X\longrightarrow\tSTP$}
	and the toric morphism
	\mbox{$\widetilde p:\widetilde X\longrightarrow\PP^A$}
	in Section~\ref{sec:toric-models}.
    
Since the Kummer modification preserves the real fans,
the cone $\tau_{\cS}$ defines an orbit closure
$\widetilde V_{\cS}\subset\tSTP$.
Set
$
\widetilde X_{\cS}\deq
\widetilde\pi^{-1}(\widetilde V_{\cS}).
$ For each $C\in\cS^{\max}$, let
$\widetilde X_{\cS,C}\subset\widetilde X_{\cS}$
be the orbit closure defined by $\sigma_{\cS,C}$.
These are the irreducible components of
$\widetilde X_{\cS}$.

For each contributing cell $C$, we will show that the map
$\widetilde X_{\cS,C}\to\PP^C$ has degree~$\mathfrak m^r g_C$, whereas the cover
$\widetilde V_{\cS}\to V_{\cS}$ has degree $\mathfrak m^r$.
The common factor cancels on descent, leaving the
contribution~$g_C$.
	
	\begin{lemma}\label{lem:degree-kummer-component-map}
		Let $C\in\cS^{\max}$ with $|C|=d+r+1$. Then
		$$
		\deg\left(
		\widetilde p^*h^{d+r}
		\cap[\widetilde X_{\cS,C}]
		\right)
		=
		\left[
		N_{\PP^A}:
		\widetilde N_X+N_{\sigma_{\cS,C}}
		\right].
		$$
	\end{lemma}
	
	\begin{proof}
		Let $\widetilde p_{\cS, C}$ be the restriction of $\widetilde p$ to
		$\widetilde X_{\cS,C}$. It factors through $\PP^C$, and its image is
		the linear subspace $\PP^C\subset\PP^A$. Since
		$|C|=d+r+1$, both $\widetilde X_{\cS,C}$ and $\PP^C$ have
		dimension $d+r$. Hence $\widetilde p_{\cS, C}$ is generically finite onto its image. The projection formula gives
		$$\deg\left(
		\widetilde p^*h^{d+r}
		\cap[\widetilde X_{\cS,C}]
		\right) = 
		\deg(\widetilde p_{\cS, C}).$$
		
		The cocharacter lattice of the dense torus of $\PP^C$ is
		$N_{\PP^A}/N_{\sigma_{\cS,C}}$, while the image of the cocharacter lattice
		$\widetilde N_X$ is
		$$
		(\widetilde N_X+N_{\sigma_{\cS,C}})/N_{\sigma_{\cS,C}}.
		$$
		Therefore
		$$
		\begin{aligned}
			\deg(\widetilde p_{\cS, C})
			&=
			\left[
			N_{\PP^A}/N_{\sigma_{\cS,C}}:
			(\widetilde N_X+N_{\sigma_{\cS,C}})/N_{\sigma_{\cS,C}}
			\right]\\
			&=
			\left[
			N_{\PP^A}:
			\widetilde N_X+N_{\sigma_{\cS,C}}
			\right].
		\end{aligned}
		$$
	\end{proof}

	Now we will compute the lattice index appearing in
	Lemma~\ref{lem:degree-kummer-component-map}.
	Recall that
	$
	\mathfrak m=g_A m$ and $	\widetilde N_X=N+\mathfrak mN_{\PP^A}
	$
	in the Kummer model described in Section~\ref{sec:toric-models}.
	\begin{lemma}\label{lem:kummer-lattice-index}
		Let $C\in\cS^{\max}$ satisfy $|C|=d+r+1$. Then
		$$
		\left[
		N_{\PP^A}:
		\widetilde N_X+N_{\sigma_{\cS,C}}
		\right]
		=
		\mathfrak m^r g_C.
		$$
	\end{lemma}
	
	\begin{proof}
		By Lemma~\ref{lem:local-effective-index},
		$$
		g_C
		=
		g_A[\Lambda_A:\Lambda_C]
		=
		g_A m_{\sigma_{\cS,C}/\tau_{\cS}}
		\mid
		g_A m
		=
		\mathfrak m.
		$$
		Since
		$$
		|N_{\Lambda_C}/N|
		=
		[M:\Lambda_C]
		=
		g_C,
		$$
		multiplication by $\mathfrak m$ annihilates
		$N_{\Lambda_C}/N$. Thus
		$
		\mathfrak mN_{\Lambda_C}\subseteq N.
		$
		
		Dualizing and quotienting the exact sequence (\ref{exactc}) by $N\subseteq N_{\Lambda_C}$,
		we obtain
		$$
		0
		\longrightarrow
		N_{\Lambda_C}/N
		\longrightarrow
		N_{\PP^C}/N
		\longrightarrow
		N_{\Sec(C)}
		\longrightarrow
		0.
		$$
		
		Since $N_{\Sec(C)}$ is free of rank
		$
		\rk N_{\Sec(C)}
		=
		|C|-d-1
		=
		r,
		$
		this exact sequence splits. Also, $\mathfrak m$ annihilates the torsion
		subgroup $N_{\Lambda_C}/N$, and thus
		$$
		\frac{N_{\PP^C}/N}
		{\mathfrak m(N_{\PP^C}/N)}
		\cong
		(N_{\Lambda_C}/N)
		\oplus
		(\ZZ/\mathfrak m\ZZ)^r.
		$$
		This group has order~$g_C\mathfrak m^r$.

		By Lemma~\ref{lem:sigma-coordinate-cone},
		$
		N_{\PP^C}
		=
		N_{\PP^A}/N_{\sigma_{\cS,C}}.
		$
		Using
		$
		\widetilde N_X
		=
		N+\mathfrak mN_{\PP^A},
		$
		we have
		$$
		\begin{aligned}
			\left[
			N_{\PP^A}:
			\widetilde N_X+N_{\sigma_{\cS,C}}
			\right]
			&=
			\left[
			N_{\PP^C}:
			N+\mathfrak mN_{\PP^C}
			\right]
			\\
			&=
			\left|
			\frac{N_{\PP^C}/N}
			{\mathfrak m(N_{\PP^C}/N)}
			\right|
			\\
			&=
			\mathfrak m^r g_C.
		\end{aligned}
		$$
	\end{proof}
	
	We can now compute the Kummer kappa class on every torus-invariant \mbox{orbit closure}.
	
	\begin{proposition}\label{prop:kummer-kappa-formula}
		For an $r$-dimensional torus-invariant orbit closure
		$\widetilde V_{\cS}\subset\tSTP$,
		$$
		\deg\left(
		\tbkp_r
		\cap[\widetilde V_{\cS}]
		\right)
		=
		\mathfrak m^r
		\sum_{\substack{
				C\in\cS^{\max}\\
				|C|=d+r+1
		}}
		g_C.
		$$
	\end{proposition}
	
	\begin{proof}
		Since $\widetilde\pi$ is flat with reduced fibers,
		$$
		[\widetilde X_{\cS}]
		=
		\sum_{C\in\cS^{\max}}
		[\widetilde X_{\cS,C}].
		$$
		Hence
		$$
		\begin{aligned}
			\deg\left(
			\tbkp_r
			\cap[\widetilde V_{\cS}]
			\right)
			&=
			\deg\left(
			\widetilde p^*h^{d+r}
			\cap[\widetilde X_{\cS}]
			\right)\\
			&=
			\sum_{C\in\cS^{\max}}
			\deg\left(
			\widetilde p^*h^{d+r}
			\cap[\widetilde X_{\cS,C}]
			\right).
		\end{aligned}
		$$
		By Corollary~\ref{cor:vanishing-nonmaximal-cells},
		only the cells satisfying
		$|C|=d+r+1$ contribute. For each such cell,
		Lemmas~\ref{lem:degree-kummer-component-map} and
		\ref{lem:kummer-lattice-index} give
		$$
		\deg\left(
		\widetilde p^*h^{d+r}
		\cap[\widetilde X_{\cS,C}]
		\right)
		=
		\mathfrak m^r g_C.
		$$
		Summing over the contributing cells proves the formula.
	\end{proof}
	
	We next descend this computation to $\STP$. The restriction
	$\mu_{\cS}: \widetilde V_{\cS} \longrightarrow V_{\cS}$
	of the Kummer lattice modification is induced by the inclusion
	$\mathfrak mN_{\Sec(A)}\subset N_{\Sec(A)}$. 
	Set
	$
	\widetilde N_{\tau_{\cS}}
	\deq
	\widetilde N_{\Sec(A)}
	\cap
	\operatorname{Span}_{\RR}(\tau_{\cS}).
	$
	Then
	$
	\widetilde N_{\tau_{\cS}}
	=
	\mathfrak mN_{\tau_{\cS}}.
	$ 
	Since $\dim V_{\cS}=r$, the induced inclusion on the cocharacter
	lattices of the dense tori of the orbit closures has index
	$$
	\begin{aligned}
		\deg(\mu_{\cS})
		&=
		\left[
		N_{\Sec(A)}/N_{\tau_{\cS}}:
		\mathfrak mN_{\Sec(A)}/\mathfrak mN_{\tau_{\cS}}
		\right]\\
		&=
		\mathfrak m^r.
	\end{aligned}
	$$
	By Proposition~\ref{prop:kummer-kappa-pullback},
	$
	\tbkp_r
	=
	\mu^*\bkp_r.
	$
	The main theorem follows immediately.
	\begin{theorem}\label{cor:intrinsic-formula}
		For every $r$-dimensional torus-invariant orbit closure
		$V_{\cS}\subset\STP$,
		$$
		\deg\left(
		\bkp_r\cap[V_{\cS}]
		\right)
		=
		\sum_{\substack{
				C\in\cS^{\max}\\
				|C|=d+r+1
		}}
		g_C.
		$$
	\end{theorem}
	
	\begin{proof}
		By the projection formula,
		$$
		\begin{aligned}
			\deg\left(
			\tbkp_r
			\cap[\widetilde V_{\cS}]
			\right)
			&=
			\deg(\mu_{\cS})
			\deg\left(
			\bkp_r
			\cap[V_{\cS}]
			\right)=
			\mathfrak m^r
			\deg\left(
			\bkp_r
			\cap[V_{\cS}]
			\right).
		\end{aligned}
		$$
		The result follows by comparing this equality with
		Proposition~\ref{prop:kummer-kappa-formula} and canceling~$\mathfrak m^r$.
	\end{proof}

	\subsection{Intersection Numbers on the Moduli Stack}
	Using the coarse descendants, it is easy to recover the intersection numbers of the kappa classes on the moduli stack $\cSTP$ from the
	coarse intersection formula obtained in
	Theorem~\ref{cor:intrinsic-formula}.
	
	Recall that $\rho:\cSTP\to\STP$ is the coarse moduli morphism.
	For an $r$-dimensional torus-invariant orbit closure
	$V_{\cS}\subset\STP$, let $\cV_{\cS}\subset\cSTP$ be the corresponding reduced
	integral substack, and let $e_{\cS}$ denote the order of the stabilizer of its geometric generic point.
	
	\begin{corollary}\label{cor:stack-intrinsic-formula}
		For every $r$-dimensional torus-invariant orbit closure
		\mbox{$V_{\cS}\subset\STP$}, we have
		$$
		\deg\left(
		\kp_r\cap[\cV_{\cS}]
		\right)
		=
		\frac{1}{e_{\cS}}
		\sum_{\substack{
				C\in\cS^{\max}\\
				|C|=d+r+1
		}}
		g_C.
		$$
	\end{corollary}
	
	\begin{proof}
		Recall that $\rho^*\bkp_r=\kp_r.$
		By the stack-to-coarse degree formula
		\cite[proof of Lemma~3.23]{AlexeevKappaBasic},
		$$
		\rho_*[\cV_{\cS}]
		=
		\frac{1}{e_{\cS}}[V_{\cS}].
		$$
		
		By the projection formula,
		$$
		\begin{aligned}
			\deg\left(
			\kp_r\cap[\cV_{\cS}]
			\right)
			&=
			\deg\left(
			\bkp_r\cap\rho_*[\cV_{\cS}]
			\right)\\
			&=
			\frac{1}{e_{\cS}}
			\deg\left(
			\bkp_r\cap[V_{\cS}]
			\right)\\
			&=
			\frac{1}{e_{\cS}}
			\sum_{\substack{
					C\in\cS^{\max}\\
					|C|=d+r+1
			}}
			g_C.
		\end{aligned}
		$$
		The last equality follows from Theorem~\ref{cor:intrinsic-formula}.
	\end{proof}

	We now recall the combinatorial description of the stabilizer
	groups of stable toric pairs from
	\cite[Sections~2.2 and~2.6]{Alexeev02}.
	We use throughout the framed, $A$-marked convention fixed above,
	so no additional symmetry group is present.
	
	Let $(\Delta_{\cS},\mathcal C_{\cS})$ be the pointed
	polyhedral complex associated with a regular subdivision $\cS$.
	For each cell $\delta\in\Delta_{\cS}$, let
	$\mathcal C_\delta$ denote its distinguished set of lattice points,
	and set
	$$
	M_\delta
	=
	M\cap\Span_{\RR}(\delta-\delta),
	\qquad
	\Lambda_\delta
	=
	\ZZ(\mathcal C_\delta-\mathcal C_\delta).
	$$
	Define
	$$
	\mathcal K_\delta
	=
	M_\delta/\Lambda_\delta,
	\qquad
	\widehat{\mathcal K}_\delta
	=
	\Hom(\mathcal K_\delta,\GG_m).
	$$
	
	By \cite[Theorem~2.6.8]{Alexeev02}, the automorphism group of the
	corresponding stable toric pair is
	$$
	\mathrm{Aut}(X,B)
	\cong
	H^0\left(
	\Delta_{\cS},
	\widehat{\mathcal K}
	[\Delta_{\cS},\mathcal C_{\cS}]
	\right)
	\cong
	\varprojlim_{\delta\in\Delta_{\cS}}
	\widehat{\mathcal K}_\delta.
	$$
	For a marked cell $C\subset A$, we write
	$\widehat{\mathcal K}_C$ for the group associated with the pointed
	cell $(\operatorname{Conv}(C),C)$ and use the same convention for
	its faces. More explicitly, the automorphism group is the subgroup of
	$\prod_{C\in\cS^{\max}}\widehat{\mathcal K}_C$
	consisting of collections $(t_C)$ satisfying
	$
	\operatorname{res}_{C,F}(t_C)
	=
	\operatorname{res}_{C',F}(t_{C'})
	$
	for every common face $F=C\cap C'$.
	Consequently, the generic stabilizer order along $\cV_{\cS}$ is
	$$
	e_{\cS}
	=
	\left|
	H^0\left(
	\Delta_{\cS},
	\widehat{\mathcal K}
	[\Delta_{\cS},\mathcal C_{\cS}]
	\right)
	\right|.
	$$
	
	Hence the stack-theoretic intersection formula,
	together with this description of the stabilizer,
	expresses the stacky intersection numbers
	entirely in terms of the combinatorial data
	of the pointed subdivision.

	\begin{corollary}\label{cor:top-kappa}
		Let $\dim\STP = n.$
		Then the normalized top kappa class satisfies
		$$
		\deg\left(
		\kp_n\cap[\cSTP]
		\right)=1.
		$$
	\end{corollary}
	
	\begin{proof}
		The trivial subdivision
		$\cS_{\mathrm{triv}}=\{A\}$
		has a single maximal cell. In the framed category,
		$$
		e_{\cS_{\mathrm{triv}}}
		=
		\left|
		\Hom(M/\Lambda_A,\GG_m)
		\right|
		=
		[M:\Lambda_A]
		=
		g_A.
		$$
		Applying Corollary~\ref{cor:stack-intrinsic-formula} gives
		$$
		\deg\left(
		\kp_n\cap[\cSTP]
		\right)
		=
		\frac{g_A}{e_{\cS_{\mathrm{triv}}}}
		=
		1.
		$$
	\end{proof}
	
	\section{Alternative Toric Computation}\label{sec:alternative-toric}
	
	We now compute the toric kappa classes using the effective
	model constructed in Section~\ref{sec:toric-models}. This alternative computation explains how replacing the ambient lattice~$M$ by the effective lattice $\Lambda_A$
	changes the resulting classes by the factor~$g_A$.
	Unlike the Kummer family, this model remains over the original base $\STP$ and need not be flat. We therefore use the toric pushforward of Minkowski
	weights rather than the flat-family construction.
	
	The computation uses the common geometry of invariant strata established in Section~\ref{subsec:common-invariant-strata}.
	
	\subsection{Toric Pushforward and Toric Kappa Classes}
	
	We recall the toric pushforward used below.
	Although a proper morphism induces a pushforward in Chow
	homology, in general it does not induce a pushforward in
	operational Chow cohomology. For the toric morphisms considered below, such a pushforward
	can be defined using Minkowski weights.
	
	Let $f:X_\Sigma\longrightarrow X_{\Sigma'}$ be a proper surjective toric morphism between complete toric varieties,
	induced by a lattice homomorphism $F:N\to N'$.
	We assume that every cone $\sigma\in\Sigma$ maps onto
	a cone of $\Sigma'$; equivalently, $f$ is equidimensional.
	Under this assumption, the pushforward of Minkowski
	weights described in \cite{GKM}, together with the
	Fulton--Sturmfels identification \cite{FS97}, defines
	$$
	f_!^{\mathrm{tor}}:
	A^p(X_\Sigma)_\QQ
	\longrightarrow
	A^{p-\dim X_\Sigma+\dim X_{\Sigma'}}
	(X_{\Sigma'})_\QQ.
	$$
	
	For $c\in A^p(X_\Sigma)_\QQ$ and
	$\tau\in\Sigma'$ satisfying
	$\dim\tau=\dim X_\Sigma-p,$ the pushforward is defined by
	$$
	\deg\left(
	f_!^{\mathrm{tor}}(c)\cap[V(\tau)]
	\right)
	=
	\sum_{\substack{
			\sigma\in\Sigma\\
			F_\RR(\sigma)=\tau\\
			\dim\sigma=\dim\tau}}
	[N_\tau:F(N_\sigma)]
	\deg\left(
	c\cap[V(\sigma)]
	\right).
	$$
	
	We apply this construction to the effective model of
	Section~\ref{sec:toric-models}. Since the effective model $\pi^{\mathrm{eff}}$ is equidimensional,
	it satisfies the assumption above. Thus the construction defines a pushforward from
	operational Chow classes on $X^{\mathrm{eff}}$ to classes
	on $\STP$.
	
	\begin{definition}
		The toric kappa class of the effective model is
		$$
		\bkpt_r
		\deq
		\pi^{\mathrm{eff, tor}}_!
		\left((p^{\mathrm{eff}*}h)^{d+r}\right)
		\in A^r(\STP)_\QQ.
		$$
	\end{definition}
	
	\subsection{The Toric Intersection Formula and Comparison}
	\begin{proposition}\label{prop:effective-toric-formula}
		For every $r$-dimensional torus-invariant orbit closure
		\mbox{$V_{\cS}\subset\STP$},
		$$
		\deg\left(
		\bkpt_r\cap[V_{\cS}]
		\right)
		=\frac{1}{g_A}
		\sum_{\substack{
				C\in\cS^{\max}\\
				|C|=d+r+1
		}}
		g_C.
		$$
	\end{proposition}
	
	\begin{proof}
		By the definition of $\bkpt_r$ and the toric
		pushforward formula,
		$$
		\deg\left(
		\bkpt_r\cap[V_{\cS}]
		\right)
		=
		\sum_{C\in\cS^{\max}}
		m_{\sigma_{\cS,C}/\tau_{\cS}}
		\deg\left(
		(p^{\mathrm{eff}})^*h^{d+r}
		\cap[X_{\cS,C}]
		\right).
		$$
		
		By Corollary~\ref{cor:vanishing-nonmaximal-cells}, only
		the cells with $|C|=d+r+1$ contribute. For such a cell,
		Lemma~\ref{lem:sigma-coordinate-cone} implies that
		$p_{\cS,C}$
		is birational onto $\PP^C$, so its contribution is $\deg((p^{\mathrm{eff}})^*h^{d+r}\cap[X_{\cS,C}])=1$.
		Moreover, by Lemma~\ref{lem:local-effective-index},
		$$
		m_{\sigma_{\cS,C}/\tau_{\cS}}
		=
		[\Lambda_A:\Lambda_C]
		=
		\frac{g_C}{g_A}.
		$$
		Therefore
		$$
		\deg\left(
		\bkpt_r\cap[V_{\cS}]
		\right)
		=
		\sum_{\substack{
				C\in\cS^{\max}\\
				|C|=d+r+1 
		}}\frac{g_C}{g_A}.
		$$
	\end{proof}
	Comparing the degrees on torus-invariant orbit closures $V_{\cS},$ we immediately get the following theorem.
	
	\begin{theorem}\label{thm:toric-comparison}
		For every $r\geq 0$,
		$$
		\bkp_r=g_A \bkpt_r
		\quad\text{in }A^r(\STP)_\QQ.
		$$
	\end{theorem}
This comparison isolates the effect of replacing the prescribed
character lattice $M$ by the effective lattice $\Lambda_A$.
	
	\section{Comparison with the GKZ Polarization}\label{sec:GKZ}
	This section identifies the first normalized kappa class $\bkp_1$ on $\STP$ with the GKZ polarization $c_1(L_{\mathrm{GKZ}})$ by comparing its degree on a torus-invariant curve with the lattice lengths of secondary edges.
	
	Fix a torus-invariant curve
	$V_{\cS}=V(\tau_{\cS})\subset\STP$,
	where $\tau_{\cS}$ is a codimension-one cone of the
	secondary fan corresponding to a regular subdivision $\cS$.
	Let $E_{\cS}\subset\Sec(A)$ be the corresponding edge,
	and let $\cT^+$ and $\cT^-$ be the regular triangulations
	corresponding to its endpoints. Note that only cells with $|C|=d+2$ contribute to the intersection formula by Lemma~\ref{lem:cell-size-bound} and Corollary~\ref{cor:vanishing-nonmaximal-cells}.
	
	\subsection{Local GKZ Differences and Lattice Indices}
	
	Let
	$C=\{c_0,\ldots,c_{d+1}\}\in\cS^{\max}$
	be a maximal cell with $|C|=d+2$, and set
	$\cT_C^\pm=\cT^\pm|_C$.
	For $j=0,\ldots,d+1$, define
	$$
	C_{\widehat j}
	\deq
	\operatorname{Conv}\bigl(C\setminus\{c_j\}\bigr).
	$$
	Then $C_{\widehat j}$ is either a $d$-simplex
	or a polytope of dimension less than $d$. We show that the index $g_C$ appearing in the intersection
	formula is the lattice length of the local GKZ difference.
	
	Consider the homogenized matrix
	$$
	\widetilde D_C
	=
	\begin{pmatrix}
		1 & 1 & \cdots & 1 \\
		c_0 & c_1 & \cdots & c_{d+1}
	\end{pmatrix}.
	$$
	For $j=0,\ldots,d+1$, let $\widetilde D_{C,\widehat j}$ be the
	matrix obtained by deleting the column corresponding to $c_j$ from $\widetilde D_C$,
	and define
	$
	m_j
	\deq
	(-1)^j\det\widetilde D_{C,\widehat j}.
	$
	Set
	$$
	m_C\deq(m_0,\ldots,m_{d+1}).
	$$
	By the cofactor identity,
	$\widetilde D_Cm_C^t=0$.
	Thus
	$$
	\sum_{j=0}^{d+1}m_j=0,
	\qquad
	\sum_{j=0}^{d+1}m_jc_j=0,
	$$
	so $m_C$ defines an affine relation among the points of $C$. Moreover,
	$$
	|m_j|
	=
	\left|\det\widetilde D_{C,\widehat j}\right|
	=
	\Vol_{\ZZ}(C_{\widehat j}).
	$$
	In particular, $m_j=0$ if and only if $C_{\widehat j}$ is
	lower-dimensional.
	
	Choose the overall sign of $m_C$ so that the restrictions of
	the fixed triangulations~$\cT^\pm$ satisfy
	$$
	(\cT_C^+)^{\max}
	=
	\{
	C_{\widehat j}:m_j>0
	\},
	\qquad
	(\cT_C^-)^{\max}
	=
	\{
	C_{\widehat j}:m_j<0
	\}.
	$$
	If $m_j=0$, then $C_{\widehat j}$ is lower-dimensional and does
	not occur as a maximal $d$-simplex in either triangulation.
	
	\begin{lemma}
		With the above sign convention,
		$$
		\varphi_{\cT_C^-}-\varphi_{\cT_C^+}
		=
		m_C.
		$$
	\end{lemma}
	
	\begin{proof}
		For $c_k\in C$, we have
		$$
		\varphi_{\cT_C^+}(c_k)
		=
		\sum_{\substack{j\ne k\\m_j>0}}
		\Vol_{\ZZ}(C_{\widehat j})
		=
		\sum_{\substack{j\ne k\\m_j>0}}m_j
		$$
		and
		$$
		\varphi_{\cT_C^-}(c_k)
		=
		\sum_{\substack{j\ne k\\m_j<0}}
		\Vol_{\ZZ}(C_{\widehat j})
		=
		-\sum_{\substack{j\ne k\\m_j<0}}m_j.
		$$
		Therefore,
		$$
		\begin{aligned}
			\bigl(
			\varphi_{\cT_C^-}
			-
			\varphi_{\cT_C^+}
			\bigr)(c_k)
			&=
			-\sum_{\substack{j\ne k\\m_j<0}}m_j
			-
			\sum_{\substack{j\ne k\\m_j>0}}m_j \\
			&=
			-\sum_{j\ne k}m_j
			=
			m_k,
		\end{aligned}
		$$
		where the last equality follows from
		$\sum_{j=0}^{d+1}m_j=0$.
	\end{proof}
	
	\begin{lemma}
		The lattice length of $m_C$ is
		$$
		\ell_{\ZZ}(m_C)=g_C.
		$$
	\end{lemma}
	
	\begin{proof}
		By definition,
		$
		\ell_{\ZZ}(m_C)
		=
		\gcd(|m_0|,\ldots,|m_{d+1}|).
		$
		Let
		$$
		\widetilde{\Lambda}_C
		=
		\ZZ\langle (1,c):c\in C\rangle
		\subset
		\ZZ\oplus M.
		$$
		By the Smith normal form, since $\widetilde D_C$
		has full row rank, the gcd of its maximal minors
		equals the index of its column lattice
		$\widetilde{\Lambda}_C$ in $\ZZ\oplus M$.
		Therefore,
		$$
		\ell_{\ZZ}(m_C)
		=
		[\ZZ\oplus M:\widetilde{\Lambda}_C].
		$$
		
		Fix $c_0\in C$. Subtracting the column $(1,c_0)$ from the other
		columns shows that $\widetilde{\Lambda}_C$ is generated by
		$(1,c_0)$ and $(0,c_j-c_0)$ for $j=1,\ldots,d+1$.
		Hence
		$$
		(\ZZ\oplus M)/\widetilde{\Lambda}_C
		\cong
		M/\ZZ(C-C).
		$$
		Therefore,
		$$
		\ell_{\ZZ}(m_C)
		=
		[M:\ZZ(C-C)]
		=
		g_C.
		$$
	\end{proof}

	\subsection{Lattice Lengths of Secondary Edges}
	
	In the previous subsection, we identified the local GKZ difference
	with $m_C$ and its lattice length with $g_C$.
	We now add these local contributions along the edge $E_{\cS}$. For a maximal cell $C$, let
	$
	\iota_C:\ZZ^C\hookrightarrow\ZZ^A
	$
	denote extension by zero.
	
	\begin{lemma}\label{Glength}
		The lattice length of $E_{\cS}$ in the lattice $L_A$ is
		$$
		\ell_{\ZZ}(E_{\cS})
		=
		\sum_{\substack{C\in\cS^{\max}\\|C|=d+2}}
		g_C.
		$$
	\end{lemma}
	
	\begin{proof}
		Let $\lambda\in L_A$ be the primitive vector in the direction from
		$\varphi_{\cT^+}$ to $\varphi_{\cT^-}$.
		By additivity of the GKZ vector over the maximal cells of $\cS$,
		$$
		\varphi_{\cT^-}-\varphi_{\cT^+}
		=
		\sum_{\substack{C\in\cS^{\max}\\|C|=d+2}}
		\iota_C(m_C).
		$$
		For every contributing cell $C$, the vector
		$\iota_C(m_C)$ lies in
		$$
		L_{A,\RR}\cap
		\bigl(\Span_{\RR}(\tau_{\cS})\bigr)^\perp.
		$$
		Since $V_{\cS}$ is a curve, this space is
		one-dimensional. Write
		$
		m_C=g_C\overline m_C,
		$
		where $\overline m_C$ is primitive.
		Since the coordinates of $\overline m_C$ have gcd $1$,
		its extension $\iota_C(\overline m_C)$ is primitive
		in $\ZZ^A$ and hence also primitive in $L_A$.
		
		Choose a height function~$\omega$ inducing $\cT^+$.
		Then $\varphi_{\cT^+}$ minimizes the pairing with~$\omega$
		on $\Sec(A)$, and $\varphi_{\cT_C^+}$ minimizes the pairing
		with $\omega|_C$ on $\Sec(C)$.
	Therefore,
		$$
		\langle\omega,\lambda\rangle>0,
		\qquad
		\langle\omega,\iota_C(m_C)\rangle>0
		$$
		for every contributing cell $C$.
		Thus $\iota_C(\overline m_C)=\lambda$.
		
		Therefore,
		$$
		\varphi_{\cT^-}-\varphi_{\cT^+}
		=
		\left(
		\sum_{\substack{C\in\cS^{\max}\\|C|=d+2}}
		g_C
		\right)\lambda.
		$$
		Since $\lambda$ is primitive, the coefficient of $\lambda$ is the
		lattice length of $E_{\cS}$.
	\end{proof}
	
	Combining the edge length formula with the main intersection formula
	identifies the operational Chow class $\bkp_1$.
	
	\begin{corollary}
		$$
		\bkp_1
		=
		c_1(L_{\mathrm{GKZ}})
		\quad\text{in }A^1(\STP)_\QQ.
		$$
	\end{corollary}
	
	\begin{proof}
		The GKZ line bundle is the ample line bundle on $\STP$
		associated with the secondary polytope $\Sec(A)$.
		For every torus-invariant curve $V_{\cS}$,
		$$
		c_1(L_{\mathrm{GKZ}})\cdot V_{\cS}
		=
		\ell_{\ZZ}(E_{\cS});
		$$
		see \cite[Proposition~6.3.8]{CLS}.
		
		On the other hand, Theorem~\ref{cor:intrinsic-formula} and
		Lemma~\ref{Glength} give
		$$
		\bkp_1\cdot V_{\cS}
		=
		\sum_{\substack{C\in\cS^{\max}\\|C|=d+2}}g_C
		=
		\ell_{\ZZ}(E_{\cS}).
		$$
		Thus the two classes have the same degree on every torus-invariant
		curve of $\STP$.
	\end{proof}	
	We illustrate how contributions from distinct maximal cells
	add along a single GKZ edge.
	\begin{example}
		\label{ex:square-center}
		
		Let~\mbox{$A=\{a_0,a_1,a_2,a_3,a_4\}\subset M \cong \ZZ^2$}, where
		$$
		\begin{aligned}
			a_0=(0,0), \;\;
			a_1=(2,0),\;\;
			a_2=(2,2),\;\;
			a_3=(0,2),\;\;
			a_4=(1,1).\;\;
		\end{aligned}
		$$
		Then $\operatorname{Conv}(A)$ is a square, and $a_4$ is an interior lattice point.
		
		Consider a regular subdivision $\cS$ having two maximal cells
		$$C_1=\{a_0,a_1,a_2,a_4\}, \quad \text{and} \quad
		C_2=\{a_0,a_2,a_3,a_4\}.
		$$
		Geometrically, the square is divided into two triangles glued along the diagonal
		$\operatorname{Conv}(a_0,a_2).$
	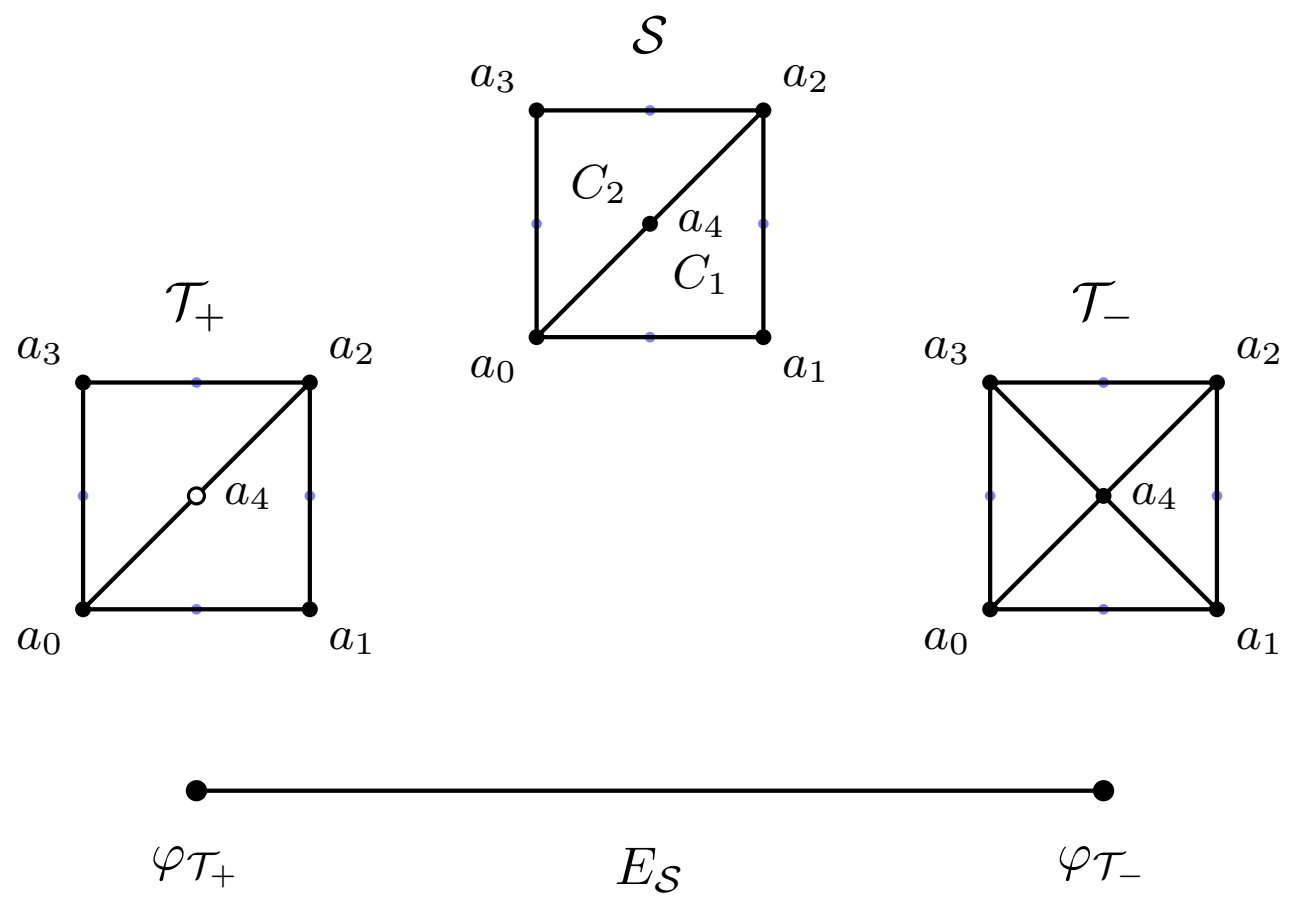
\begin{figure}[htbp]
		\centering
		\begin{tikzpicture}[
			x=1cm,y=1cm,
			line cap=round,
			line join=round
			]
			\node at (0,3.8) {$\mathcal S$};
			\node at (-3,2.0) {$\mathcal T_+$};
			\node at ( 3,2.0) {$\mathcal T_-$};
			
			\foreach \xshift/\yshift/\kind in {
				-3.75/0/0,-0.75/1.8/1,2.25/0/2
			}{
				\begin{scope}[
					shift={(\xshift,\yshift)},
					x=0.75cm,y=0.75cm,
					font=\small
					]
					\coordinate (a0) at (0,0);
					\coordinate (a1) at (2,0);
					\coordinate (a2) at (2,2);
					\coordinate (a3) at (0,2);
					\coordinate (a4) at (1,1);
					
					\foreach \i in {0,1,2}{
						\foreach \j in {0,1,2}{
							\fill[blue!50] (\i,\j) circle (1pt);
						}
					}
					
					\draw[thick] (a0)--(a1)--(a2)--(a3)--cycle;
					\draw[thick] (a0)--(a2);
					
					\ifnum\kind=2\relax
					\draw[thick] (a1)--(a3);
					\fi
					
					\ifnum\kind=0\relax
					\draw[fill=white,line width=0.6pt]
					(a4) circle (1.5pt);
					\else
					\fill (a4) circle (1.5pt);
					\fi
					
					\fill (a0) circle (1.5pt) node[below left] {$a_0$};
					\fill (a1) circle (1.5pt) node[below right] {$a_1$};
					\fill (a2) circle (1.5pt) node[above right] {$a_2$};
					\fill (a3) circle (1.5pt) node[above left] {$a_3$};
					\node[right=5pt,inner sep=0pt]
					at (a4) {$a_4$};
					
					\ifnum\kind=1\relax
					\node at (1.45,0.55) {$C_1$};
					\node at (0.55,1.35) {$C_2$};
					\fi
				\end{scope}
			}
			
			\draw[thick] (-3,-1.2)--(3,-1.2)
			node[midway,below=7pt]
			{$E_{\mathcal S}$};
			
			\fill (-3,-1.2) circle (2pt)
			node[below=7pt] {$\varphi_{\mathcal T_+}$};
			\fill ( 3,-1.2) circle (2pt)
			node[below=7pt] {$\varphi_{\mathcal T_-}$};
		\end{tikzpicture}
		\caption{
			The subdivision $\mathcal S$ corresponds to the edge
			$E_{\mathcal S}$, whose endpoints are the GKZ vectors of
			$\mathcal T_+$ and $\mathcal T_-$.
			The open circle indicates that $a_4$ is not used as a vertex of $\mathcal T_+$.
		}
		\label{fig:square-center-edge}
	\end{figure}
	
Both cells have the affine relation $a_0+a_2=2a_4$.
		Moreover,
		$$
		\Lambda_{C_1}=\Lambda_{C_2}
		=\{(p,q)\in\mathbb Z^2:p - q \in 2\ZZ\}.
		$$
		Hence
		$
		g_{C_1}=g_{C_2}=2.
		$
		
		Each local secondary polytope is a segment.
		After extension by zero to $\mathbb R^A$,
		both segments are parallel to
		$$
		\lambda=(-1,0,-1,0,2)\in L_A.
		$$
		By the face decomposition of the secondary
		polytope, the corresponding face is an edge $E_{\cS}$.
		Thus $V_{\cS}$ is a curve.
		
		Let $\cT_+$ be the triangulation along
		the diagonal $a_0a_2$, omitting $a_4$,
		and let $\cT_-$ be the triangulation
		consisting of the four triangles meeting
		at $a_4$.
		
		Their GKZ vectors are
		$$
		\varphi_{\cT_+}=(8,4,8,4,0),
		\qquad
		\varphi_{\cT_-}=(4,4,4,4,8).
		$$
		Each maximal cell contributes $2\lambda$
		to their difference, so
		$$
		\varphi_{\cT_-}-\varphi_{\cT_+}
		=2\lambda+2\lambda=4\lambda.
		$$
		Since $\lambda$ is primitive in $L_A$,
		the corresponding edge has lattice length
		$$
		\ell_{\mathbb Z}(E_{\cS})=4.
		$$
		
		On the other hand, our intersection formula
		gives
		$$
		\deg\left(
		\bkp_1\cap[V_{\cS}]
		\right)
		=g_{C_1}+g_{C_2}
		=2+2=4.
		$$
		Thus the two local contributions recover
		the lattice length of the GKZ edge.
	\end{example}

\section{\texorpdfstring
	{The Interval Configuration $A_N=\{0,\ldots,N\}$}
	{The Interval Configuration}}\label{sec:interval}
	
	The main intersection formula determines the classes through
	their degrees on orbit closures. For interval configurations, we turn these degrees into explicit formulas in the Chow ring and compute all top intersection products.
	Let~\mbox{$A_N=\{0,1,\ldots,N\}\subset\mathbb Z$}.
	The secondary fan is simplicial, so the associated toric variety has a
	Stanley--Reisner presentation for its rational Chow ring. We use this
	presentation to derive closed formulas for the normalized kappa classes. The case $N=2$ was considered in
	Example~\ref{ex:three-points}.
	We assume $N\geq3$ below, so that the
	displayed ray generators are primitive.

	\subsection{The Secondary Fan and Its Chow Ring}
	\label{subsec:interval-chow}
	
	Let \mbox{$A_N=\{0,1,\ldots,N\}\subset\mathbb Z$} for $N\geq3,$ and let
	$\STP_N=X_{\operatorname{Sec}(A_N)}$ be the toric variety associated
	with its secondary polytope. We write
	$n=N-1=\dim\STP_N$ and $[n]=\{1,\ldots,n\}.$
	
	By \cite[Chapter~7, Section~3A]{GKZ}, the secondary polytope of
	$A_N$ is combinatorially equivalent to an $n$-dimensional cube.
	Its faces can be described by pairs of disjoint subsets
	$U,B\subseteq[n]$, where $U$ is the set of omitted interior points and $B$ is the set of breakpoints. Points outside $U\cup B$ remain marked within a cell
	but are not prescribed as breakpoints. We denote the corresponding marked subdivision by~$\mathcal S(U,B)$.
	
	More precisely, if $B=\{b_1<\cdots<b_k\}$,
	set $b_0=0$ and $b_{k+1}=N$. Then
	$\mathcal S(U,B)$ consists of the marked cells
	$$
	\left(
	[b_{\ell-1},b_\ell],
	(A_N\cap[b_{\ell-1},b_\ell])\setminus U
	\right),
	\qquad 1\leq\ell\leq k+1.
	$$
	
	Requiring an interior point to be omitted or to be a
	breakpoint defines a facet of the secondary polytope. For each $j\in[n],$ let $E_j$ and $F_j$ denote the rays corresponding
	respectively to omitting $j$ and prescribing a breakpoint at $j,$ and
	write $v_{E_j}$ and $v_{F_j}$ for their primitive generators in the
	lattice of the secondary fan. Thus
	$$
	\sigma_{U,B}
	=
	\operatorname{Cone}
	\bigl(
	\{v_{E_j}:j\in U\}
	\cup
	\{v_{F_j}:j\in B\}
	\bigr),
	\qquad
	U\cap B=\varnothing,
	$$
	and the corresponding torus-invariant orbit closure has dimension
	$$
	\dim V_{\mathcal S(U,B)}
	=
	n-|U|-|B|.
	$$
	
	To compute the Chow ring, we first record the primitive ray generators
	in convenient lattice coordinates. As in Section~\ref{subsec:kummer-model}, the lattice of
	the secondary fan is naturally $L_{A_N}^{\vee}$, where
	$$
	L_{A_N}
	=
	\left\{
	(c_0,\ldots,c_N)\in\mathbb Z^{N+1}
	:
	\sum_{i=0}^N c_i
	=
	\sum_{i=0}^N ic_i
	=
	0
	\right\}.
	$$
	Let $\epsilon_0,\ldots,\epsilon_N$ denote the standard basis of
	$\mathbb Z^{N+1}.$ The elementary relations
	$$
	b_j
	=
	\epsilon_{j-1}-2\epsilon_j+\epsilon_{j+1},
	\qquad
	1\leq j\leq n,
	$$
	form an integral basis of $L_{A_N}.$
	Let $u_1,\ldots,u_n$ be the dual basis of $L_{A_N}^\vee,$ with the
	convention $u_0=u_N=0.$
	
	With respect to this basis, the facet equations
	\cite[Chapter~7, (3.4)--(3.5)]{GKZ} give
	$$
	v_{F_j}=u_j,
	\qquad
	v_{E_j}=u_{j-1}-2u_j+u_{j+1},
	\qquad
	1\leq j\leq n.
	$$
	For $N\geq3,$ these are primitive lattice vectors.
	
	We can now describe the rational Chow ring of $\STP_N$.
	
	\begin{proposition}
		\label{prop:interval-chow-ring}
		Let $x_j$ and $y_j$ denote the classes of the torus-invariant
		divisors corresponding to $E_j$ and $F_j,$ respectively. Then
		$$
		A^*(\STP_N)_{\mathbb Q}
		\simeq
		\frac{
			\mathbb Q[x_1,\ldots,x_n, y_1, \ldots, y_n]
		}{
			\left(x_jy_j, 
			y_j -2x_j+x_{j-1}+x_{j+1}
			\right)_{j=1}^{n}
		},
		$$
		where $x_0=x_N=0.$
	\end{proposition}
	
	\begin{proof}
		The opposite facets of the cube correspond to the pairs
		$E_j,F_j,$ so the Stanley--Reisner ideal is
		$(x_1y_1,\ldots,x_ny_n).$
		The ray generators give the linear relations
		$$
		y_j-2x_j+x_{j-1}+x_{j+1}=0,
		\qquad
		1\leq j\leq n.
		$$
		Eliminating the variables $y_j$ from the Stanley--Reisner
		relations gives
		$$
		2x_j^2=x_j(x_{j-1}+x_{j+1}),
		\qquad
		1\leq j\leq n,
		$$
		with $x_0=x_N=0.$
	\end{proof}
	
	This presentation identifies $A^*(\STP_N)_{\mathbb Q}$ with the
	specialization at $q=1$ of the finite Klyachko algebra $K_N$ of
	\cite{NadeauTewari}. Hence the squarefree monomials
	\mbox{$x_I:=\prod_{i\in I}x_i$}, for $I\subseteq[n],$ form a basis of
	$A^*(\STP_N)_{\mathbb Q}.$
	For $J\subseteq[n],$ we similarly write
	\mbox{$y_J:=\prod_{j\in J}y_j$}.
	
	\begin{lemma}
		\label{lem:interval-diagonal-pairing}
		Let $I,J\subseteq[n]$ with $|I|=|J|,$ and write
		$I=B_1\sqcup\cdots\sqcup B_t$ for the decomposition of $I$
		into maximal consecutive blocks. Then
		$$
		\int_{\STP_N}x_Iy_{J^c}
		=
		\begin{cases}
			\displaystyle
			\prod_{a=1}^{t}\frac{1}{|B_a|+1},
			& I=J,\\[8pt]
			0,
			& I\neq J,
		\end{cases}
		$$
		where $J^c=[n]\setminus J.$
	\end{lemma}
	
	\begin{proof}
		If $I\neq J,$ then $I\cap J^c\neq\varnothing,$ so
		$x_Iy_{J^c}=0$ by the Stanley--Reisner relations.
		
		Suppose $I=J.$ The product $x_Iy_{I^c}$ corresponds to the
		maximal cone $\sigma_{I,I^c}.$ With respect to the basis
		$u_1,\ldots,u_n,$ the determinant defining its multiplicity
		decomposes according to the maximal consecutive blocks
		$B_1,\ldots,B_t$ of $I.$ A block of length $m$ contributes the
		type-$A_m$ tridiagonal matrix
		$$
		\begin{pmatrix}
			-2&1&&\\
			1&-2&\ddots&\\
			&\ddots&\ddots&1\\
			&&1&-2
		\end{pmatrix},
		$$
		whose determinant has absolute value $m+1.$ Thus
		$$
		\operatorname{mult}(\sigma_{I,I^c})
		=
		\prod_{a=1}^t(|B_a|+1).
		$$
		The toric intersection formula gives the result.
	\end{proof}
	
	\subsection{Explicit Formulas for Kappa Classes}
	\label{subsec:interval-kappa}
	
	We now apply our main intersection formula to determine the normalized
	kappa classes explicitly. We first express them in the squarefree basis
	and then rewrite the result in terms of power sums.
	
	\begin{proposition}
		\label{prop:interval-explicit-kappa}
		For $1\leq r\leq n,$
		$$
		\bkp_r
		=
		(r+1)\sum_{j=1}^{n-r+1}
		x_jx_{j+1}\cdots x_{j+r-1}.
		$$
		Moreover, $\bkp_0=N.$
	\end{proposition}
	
	\begin{proof}
		Write
		$$
		\bkp_r
		=
		\sum_{\substack{I\subseteq[n]\\|I|=r}}
		a_Ix_I
		$$
		in the squarefree basis. By
		Lemma~\ref{lem:interval-diagonal-pairing},
		$$
		\int_{\STP_N}\bkp_r\,y_{I^c}
		=
		a_I\int_{\STP_N}x_Iy_{I^c}.
		$$
		
		The cone corresponding to
		$\mathcal S(\varnothing,I^c)$ is generated by
		$\{v_{F_j}:j\in I^c\}.$ Since \mbox{$v_{F_j}=u_j$}, this cone is
		unimodular, and hence
		$$
		y_{I^c}\cap[\STP_N]
		=
		[V_{\mathcal S(\varnothing,I^c)}].
		$$
		In the subdivision $\mathcal S(\varnothing,I^c),$ the maximal
		cells containing at least three marked points are indexed by the
		maximal consecutive blocks $B$ of $I.$ The cell $C$ corresponding
		to $B$ contains $|B|+2$ consecutive marked points and has
		$g_C=1.$ The remaining maximal cells are unit intervals.
		
		Since $|I|=r,$ the main intersection formula gives
		$$
		\int_{\STP_N}\bkp_r\,y_{I^c}
		=
		\begin{cases}
			1,& I\text{ is consecutive},\\
			0,& \text{otherwise}.
		\end{cases}
		$$
		Indeed, a maximal cell contributes precisely when it contains
		$r+2$ marked points, which occurs exactly when $I$ consists of
		a single consecutive block.
		
		If $I$ is consecutive, then
		$$
		\int_{\STP_N}x_Iy_{I^c}
		=
		\frac{1}{r+1}
		$$
		by Lemma~\ref{lem:interval-diagonal-pairing}.
		Hence $a_I=r+1$ for consecutive $I,$ and $a_I=0$ otherwise.
		
		Finally, $\bkp_0$ equals the normalized volume of $[0,N],$
		which is $N.$
	\end{proof}
	
	\begin{corollary}
		\label{cor:interval-kappa-power-sum}
		For $1\leq r\leq n,$
		$$
		\bkp_r
		=
		(r+1)\sum_{i=1}^{n}x_i^r.
		$$
	\end{corollary}
	
	\begin{proof}
		Iterating the $q=1$ specialization of
		\cite[Lemma~3.6]{NadeauTewari} and summing over the generators
		in the finite Klyachko algebra gives
		$$
		\sum_{i=1}^{n}x_i^r
		=
		\sum_{j=1}^{n-r+1}
		x_jx_{j+1}\cdots x_{j+r-1}.
		$$
		The result follows from
		Proposition~\ref{prop:interval-explicit-kappa}.
	\end{proof}
	
	\subsection{Product Formula}
	\label{subsec:interval-products}
	
	Set $p_r=\sum_{i=1}^{n}x_i^r.$
	By Corollary~\ref{cor:interval-kappa-power-sum}, it suffices to
	compute top products of the power sums. We use the cyclic sum rule
	for the Klyachko algebra \cite[Proposition~5.3]{NadeauTewari}.
	
	\begin{lemma}
		\label{lem:interval-symmetric-degree}
		Let $F\in\mathbb Q[z_1,\ldots,z_N]$ be a symmetric homogeneous
		polynomial of degree $n=N-1.$ Then
		$$
		\int_{\STP_N}
		F(x_1,\ldots,x_n,0)
		=
		\frac{F(1,\ldots,1)}{N^2}.
		$$
	\end{lemma}
	
	\begin{proof}
		Set $T=x_1\cdots x_n.$ By
		Lemma~\ref{lem:interval-diagonal-pairing},
		$$
		\int_{\STP_N}T=\frac1N.
		$$
		For $G\in\mathbb Q[x_1,\ldots,x_n]$ homogeneous of degree $n,$
		let $\operatorname{Top}(G)$ denote the coefficient of $T$ in its
		reduction in the squarefree basis.
		
		Expand $F$ into monomials and group their $N$-tuples of exponents
		into cyclic orbits. Since the total degree is $N-1,$ every exponent
		vector has a zero coordinate. Moreover, every cyclic orbit has
		exactly $N$ elements, since any nontrivial repetition factor would
		divide both $N$ and $N-1.$
		
		Since $F$ is symmetric, all monomials in a fixed cyclic orbit have
		the same coefficient. The coefficient formula and cyclic sum rule in
		\cite[Theorem~4.11 and Proposition~5.3]{NadeauTewari},
		specialized to $q=1$, give
		$$
		\sum_{\substack{
				\beta\text{ in the cyclic orbit}\\
				\beta_N=0}}
		\operatorname{Top}
		\left(
		x_1^{\beta_1}\cdots x_n^{\beta_n}
		\right)
		=
		1.
		$$
		Hence each orbit contributes its common coefficient to
		$\operatorname{Top}(F(x_1,\ldots,x_n,0)),$ while it contributes
		$N$ times that coefficient to $F(1,\ldots,1).$ Therefore
		$$
		\operatorname{Top}
		\left(F(x_1,\ldots,x_n,0)\right)
		=
		\frac{F(1,\ldots,1)}{N}.
		$$
		Multiplying by $\int_{\STP_N}T=1/N$ gives the result.
	\end{proof}
	
	\begin{theorem}
		\label{thm:interval-product-formula}
		Let $r_1,\ldots,r_s$ be positive integers satisfying
		$r_1+\cdots+r_s=N-1.$ Then
		$$
		\int_{\STP_N}
		\prod_{i=1}^{s}\bkp_{r_i}
		=
		N^{s-2}
		\prod_{j=1}^{s}(r_j+1).
		$$
	\end{theorem}
	
	\begin{proof}
		By Corollary~\ref{cor:interval-kappa-power-sum},
		$$
		\int_{\STP_N}
		\prod_{j=1}^{s}\bkp_{r_j}
		=
		\left(
		\prod_{j=1}^{s}(r_j+1)
		\right)
		\int_{\STP_N}
		\prod_{j=1}^{s}p_{r_j}.
		$$
		Apply Lemma~\ref{lem:interval-symmetric-degree} to
		$$
		F(z_1,\ldots,z_N)
		=
		\prod_{j=1}^{s}
		\left(
		\sum_{i=1}^{N}z_i^{r_j}
		\right).
		$$
		This is symmetric and homogeneous of degree
		$r_1+\cdots+r_s=N-1,$ and
		$$
		F(x_1,\ldots,x_n,0)
		=
		\prod_{j=1}^{s}p_{r_j},
		\qquad
		F(1,\ldots,1)=N^s.
		$$
		Hence
		$$
		\int_{\STP_N}
		\prod_{j=1}^{s}p_{r_j}
		=
		N^{s-2},
		$$
		which proves the formula.
	\end{proof}

	\bibliographystyle{alpha}

\begin{thebibliography}{GKM09}
		
		\bibitem[AB06]{AlexeevBrion06}
		Valery Alexeev and Michel Brion.
		\newblock Stable spherical varieties and their moduli.
		\newblock {\em International Mathematics Research Papers}, 2006:1--57, 2006.
		
		\bibitem[Ale02]{Alexeev02}
		Valery Alexeev.
		\newblock Complete moduli in the presence of semiabelian group action.
		\newblock {\em Annals of Mathematics. Second Series}, 155(3):611--708, 2002.
		
		\bibitem[Ale06]{Alexeev06}
		Valery Alexeev.
		\newblock Higher-dimensional analogues of stable curves.
		\newblock In {\em Proceedings of the International Congress of Mathematicians,
			Madrid 2006, Volume II}, pages 515--536. European Mathematical Society,
		Z{\"u}rich, 2006.
		
		\bibitem[Ale15]{Weighted}
		Valery Alexeev.
		\newblock {\em Moduli of weighted hyperplane arrangements}.
		\newblock Advanced Courses in Mathematics. CRM Barcelona.
		Birkh\"auser/Springer, Basel, 2015.
		
		\bibitem[Ale25]{Ale23}
		Valery Alexeev.
		\newblock Kappa classes on {KSBA} spaces.
		\newblock {\em Moduli}, 2:Paper No. e4, 14, 2025.
		
		\bibitem[Ale26]{AlexeevKappaBasic}
		Valery Alexeev.
		\newblock Tautological classes on {KSBA} spaces, 2026.
		\newblock arXiv:2607.19251.
		
		\bibitem[AM16]{AscherMolcho16}
		Kenneth Ascher and Samouil Molcho.
		\newblock Logarithmic stable toric varieties and their moduli.
		\newblock {\em Algebraic Geometry}, 3(3):296--319, 2016.
		
		\bibitem[AP26]{AlexeevPardini}
		Valery Alexeev and Rita Pardini.
		\newblock Explicit compactifications of moduli spaces of {Campedelli} and
		{Burniat} surfaces.
		\newblock {\em Annali della Scuola Normale Superiore di Pisa, Classe di
			Scienze}, 27(1):445--490, 2026.
		
		\bibitem[BS92]{BilleraSturmfels92}
		Louis~J. Billera and Bernd Sturmfels.
		\newblock Fiber polytopes.
		\newblock {\em Annals of Mathematics}, 135(3):527--549, 1992.
		
		\bibitem[CLS11]{CLS}
		David~A. Cox, John~B. Little, and Henry~K. Schenck.
		\newblock {\em Toric Varieties}, volume 124 of {\em Graduate Studies in
			Mathematics}.
		\newblock American Mathematical Society, Providence, RI, 2011.
		
		\bibitem[Don20]{Don23}
		Simon Donaldson.
		\newblock Fredholm topology and enumerative geometry: reflections on some words
		of {M}ichael {A}tiyah.
		\newblock In {\em Proceedings of the {G}\"okova {G}eometry-{T}opology
			{C}onferences 2018/2019}, pages 1--31. Int. Press, Somerville, MA, 2020.
		
		\bibitem[FM81]{FM}
		William Fulton and Robert MacPherson.
		\newblock Categorical framework for the study of singular spaces.
		\newblock {\em Mem. Amer. Math. Soc.}, 31(243):vi+165, 1981.
		
		\bibitem[FS97]{FS97}
		William Fulton and Bernd Sturmfels.
		\newblock Intersection theory on toric varieties.
		\newblock {\em Topology}, 36(2):335--353, 1997.
		
		\bibitem[Ful98]{Ful}
		William Fulton.
		\newblock {\em Intersection Theory}, volume~2 of {\em Ergebnisse der Mathematik
			und ihrer Grenzgebiete. 3. Folge}.
		\newblock Springer-Verlag, Berlin, second edition, 1998.
		
		\bibitem[GKM09]{GKM}
		Andreas Gathmann, Michael Kerber, and Hannah Markwig.
		\newblock Tropical fans and the moduli spaces of tropical curves.
		\newblock {\em Compos. Math.}, 145(1):173--195, 2009.
		
		\bibitem[GKZ94]{GKZ}
		I.~M. Gel'fand, M.~M. Kapranov, and A.~V. Zelevinsky.
		\newblock {\em Discriminants, resultants, and multidimensional determinants}.
		\newblock Mathematics: Theory \& Applications. Birkh\"auser Boston, Inc.,
		Boston, MA, 1994.
		
		\bibitem[Jia22]{JiangVFC}
		Yunfeng Jiang.
		\newblock The virtual fundamental class for the moduli space of surfaces of
		general type, 2022.
		\newblock arXiv:2206.00575.
		
		\bibitem[Jia26]{JiangEnumerativeKSBA}
		Yunfeng Jiang.
		\newblock Enumerative geometry on {KSBA} moduli spaces, 2026.
		\newblock arXiv:2605.19923.
		
		\bibitem[Kol23]{Kollar23}
		J{\'a}nos Koll{\'a}r.
		\newblock {\em Families of Varieties of General Type}, volume 231 of {\em
			Cambridge Tracts in Mathematics}.
		\newblock Cambridge University Press, 2023.
		\newblock With the collaboration of Klaus Altmann and S{\'a}ndor J. Kov{\'a}cs.
		
		\bibitem[Laf03]{Lafforgue2003}
		Laurent Lafforgue.
		\newblock {\em Chirurgie des grassmanniennes}, volume~19 of {\em CRM Monograph
			Series}.
		\newblock American Mathematical Society, Providence, RI, 2003.
		
		\bibitem[Mil86]{Mil}
		Edward~Y. Miller.
		\newblock The homology of the mapping class group.
		\newblock {\em Journal of Differential Geometry}, 24(1):1--14, 1986.
		
		\bibitem[Mol21]{Molcho}
		Sam Molcho.
		\newblock Universal stacky semistable reduction.
		\newblock {\em Israel Journal of Mathematics}, 242(1):55--82, 2021.
		
		\bibitem[Mor87]{Mor}
		Shigeyuki Morita.
		\newblock Characteristic classes of surface bundles.
		\newblock {\em Invent. Math.}, 90(3):551--577, 1987.
		
		\bibitem[Mum83]{Mum}
		David Mumford.
		\newblock Towards an enumerative geometry of the moduli space of curves.
		\newblock In {\em Arithmetic and geometry, {V}ol. {II}}, volume~36 of {\em
			Progr. Math.}, pages 271--328. Birkh\"auser Boston, Boston, MA, 1983.
		
		\bibitem[MW07]{MW07}
		Ib~Madsen and Michael Weiss.
		\newblock The stable moduli space of {R}iemann surfaces: {M}umford's
		conjecture.
		\newblock {\em Ann. of Math. (2)}, 165(3):843--941, 2007.
		
		\bibitem[NT25]{NadeauTewari}
		Philippe Nadeau and Vasu Tewari.
		\newblock A {$q$}-deformation of an algebra of {Klyachko} and {Macdonald}'s
		reduced word formula.
		\newblock {\em Transactions of the American Mathematical Society},
		378(12):8821--8870, 2025.
		
		\bibitem[Vis89]{Vis}
		Angelo Vistoli.
		\newblock Intersection theory on algebraic stacks and on their moduli spaces.
		\newblock {\em Invent. Math.}, 97(3):613--670, 1989.
		
	\end{thebibliography}

\end{document}